\magnification=\magstep1
\input epsf

~~\vskip 1in

\centerline{\bf Cell Decomposition and Compactification}

\centerline{\bf of}

\centerline {\bf Riemann's Moduli Space }

\centerline{\bf in}

\centerline{\bf Decorated Teichm\"uller Theory}

\vskip .2in

\centerline {\bf R. C. Penner}

\centerline{Departments of Mathematics and Physics/Astronomy}

\centerline{University of Southern California}

\centerline {Los Angeles, CA 90089}

\vskip .2in

\vskip .2in

\centerline{May 22, 2003}

\vskip .2in

\leftskip .8in\rightskip .8in

\noindent {\bf Abstract}~~This survey covers earlier work of the author as well as recent work on Riemann's moduli space, its canonical
cell decomposition and compactification, and the related operadic structure of arc complexes.

\leftskip=0ex\rightskip=0ex

~~\vskip .2in

\centerline{ \bf Introduction}\vskip .2in

\noindent Riemann's moduli space $M=M(F)$ of a surface $F$ has an essentially canonical cell decomposition ([Ha]-[Ko]-[St] or 
[P1]) and
admits various interesting compactifications, some of which are less or more compatible with the cell decomposition.  There are
furthermore several related operads based upon supersets of $M$, that is, there is a natural composition of appropriate (possibly
degenerate) surfaces.  These
structures on moduli space are useful tools in several guises on the interfaces of mathematics and physics.  For instance and
essentially by definition, a conformal field theory leads to an algebra over the operad of uncompactified moduli space, where the
operad structure is induced by gluing surfaces along boundary components, and algebras over this operad give Batalin-Vilkovisky
structures [Ge]. In practice, one is furthermore interested in
compactifications, or partial compactifications, of $M$; for instance, the Deligne-Mumford compactification [DM] and its operad in an
appropriate setting correspond to Gromov-Witten theory and quantum cohomology.  These remarks emphasize the importance of understanding
compactifications of $M$.

\vskip .1in

\noindent A basic fact about moduli space $M=M(F)$ is that it admits a natural cell decomposition; there are several incarnations of
cell decompositions depending upon the type of surface $F$, and we distinguish between ``punctured '' (surfaces without
boundary and with at least one puncture) and ``bordered'' (surfaces with non-empty boundary which may also have punctures).
As is well-known in the punctured case, $M$ is a non-compact orbifold, which comes equipped with a canonical cell decomposition, where
the cells in the decomposition are in one-to-one correspondence with suitable classes of graphs, or equivalently, with suitable classes
of arc families in $F$; this cell decomposition does not extend to the DM compactification in any known way.    

\vskip .1in

\noindent In the setting of this paper in the case of a punctured surface $F_1$, upon choosing a
distinguished puncture of $F_1$, there is a canonical compactification $Arc(F_1)$ of $M(F_1)$, which arises as the quotient of
the simplicial completion of the natural cell decomposition.  In the setting of this paper in the case of a bordered surface $F_2$,
there is a natural action of the positive reals ${\bf R}_{>0}$ on $M(F_2)$ which is generated by scaling all of the hyperbolic lengths
of the boundary geodesics, and we shall let $M(F_2)/{\bf R}_{>0}$ denote the quotient.
There is a canonical compactification
$Arc(F_2)$ of a space proper homotopy equivalent to
$M(F_2)/{\bf R}_{>0}$, where $Arc(F_2)$ comes equipped with a canonical cell decomposition. 
We had previously conjectured [P3] that in the case of a bordered surface $F_2$, $Arc(F_2)$ is PL-homeomorphic to a sphere and shall
discuss in $\S$5 this ``sphericity conjecture''.  This generalizes to arbitrary bordered surfaces a classical
combinatorial fact about triangulations of polygons as we shall see.  In fact, this sphericity of
$Arc(F_2)$ for any bordered surface $F_2$, in turn, implies that the compactification $Arc(F_1)$ of moduli space is an orbifold, for any
punctured surface $F_1$ as we shall also see.  In particular, $Arc(F_1)$ is conjecturally a new orbifold compactification of $M(F_1)$
for a punctured surface $F_1$, where the cell decomposition of $M(F_1)$ extends to $Arc(F_1)$.  

\vskip .1in

\noindent Other aspects of the geometry and
combinatorics of these new compactifications had already been described and will be recalled here as well.
For instance: whereas the DM compactification is regarded in  our formalism as ``forgetting'' the geometric structure on sub-annuli
of $F$, our compactification allows one to ``forget'' geometric structure on {\sl any} sub-surface embedded in $F$; a point in our
compactification gives a canonical decomposition of $F$ into respective geometric and topological sub-surfaces, and the
topological part is non-empty if and only if the point is ideal; whereas
Thurston's compactification (see [PH] for instance) of Teichm\"uller space records only those hyperbolic lengths which diverge, our
compactification records only those hyperbolic lengths which do not diverge.  

\vskip .1in

\noindent On the algebraic side, Sullivan [CS], [Su], [CJ], [Vo] has introduced a ``string prop'', which is intended to formalize in the
language of traditional algebraic topology certain algebraic aspects of physical and mathematical theories involving colliding strings.
This prop is intended to be ``universal'' in the sense that it should act on the loop space of any manifold.  The precise definition of
the appropriate compactification is apparently still being formulated.

\vskip .1in

\noindent In our
setting, there are natural operad structures on arc complexes in the bordered case which are closely connected to the string prop.
These  ``arc operads'' have been studied in [KLP] (and are recalled here in $\S$7 with a more thorough survey given in Kaufmann's
paper [K2] in this volume, which also surveys his independent work from [K1]). Thus, this paper both surveys older results [P1-P3] (in
$\S\S$1,2,6) and sketches (and even excerpts at points) the newer works [P4] on bordered surfaces (in
$\S$3),  on sphericity (in
$\S$5), and [KLP] on the operad structure (in $\S$7).

\vskip .1in

\noindent We have attempted to include enough background material to provide a meaningful survey, for instance, including many of the
principal identities of decorated Teichm\"uller theory in $\S$2 and sketches of proofs, where appropriate.  The intention has been to
emphasize the background material for the compactification and for the operad structure, and there are other aspects [P5] of decorated
Teichm\"uller theory (notably matrix models, universal constructions, homeomorphisms of the circle and Lie algebras, wavelets and
Fourier transform) which are not discussed here.

\vskip .1in

\noindent This paper is organized as follows.  $\S$1 gives basic definitions and compares the conformal treatment of the cell
decomposition of moduli space in the punctured case with the hyperbolic treatment of the cell decomposition.   $\S$2 describes the
several coordinates systems of decorated Teichm\"uller theory and develops from first principles several deep results.
 As mentioned before,
$\S\S$3-4 survey our newer work on bordered surfaces, and $\S$5 gives all the details of a proof of the sphericity conjecture for
punctured polygons.
$\S$6 contains older work that explains the degeneration of
hyperbolic structure in our new compactification, and $\S$7 surveys the arc operad and its relationship with the
sphericity conjecture.  A final appendix is intended to give a mathematical overview of folding problems
in computational biology, employs the arc complex constructions to give new combinatorial models
for certain biopolymers, and ends with a speculative remark.

\vskip .3in

\noindent {\bf 1. Definitions and cell decomposition for punctured surfaces}\vskip .2in

\noindent Consider a smooth surface $F=F_{g,r}^s$ of genus $g\geq 0$ with $r\geq 0$ labeled boundary components
and
$s\geq 0$ labeled punctures, where 
$6g-7+4r+2s\geq 0$.  The surface $F$ is said to be {\it bordered} if
$r\neq 0$.  The {\it pure mapping class group} $PMC=PMC(F)$
of
$F$ is the group of all isotopy classes of orientation-preserving homeomorphisms of $F$ which pointwise fix each boundary component
and each puncture, where the isotopy is likewise required to pointwise fix these sets. 

\vskip .1in

\noindent In the unbordered case $r=0$, define the classical 
{\it Teichm\"uller space} $T_{g}^s$ of $F$ to be the
space of all complete finite-area metrics of constant Gauss
curvature ${-1}$ (``hyperbolic metrics'') on $F$ 
modulo push-forward by diffeomorphisms fixing each puncture which are isotopic to
the identity relative to the punctures.   $T_{g}^s$ is homeomorphic to an open ball of real dimension
$6g-6+2s$, the mapping class group $PMC$ acts on $T_{g}^s$ by push-forward of metric under
representative diffeomorphism, and the
quotient $M_{g}^s=T_{g}^s/PMC$ is by definition {\it Riemann's
moduli space} of $F$.

\vskip .1in

\noindent In the unbordered punctured case, where $r=0$ and $s>0$, there is a  
trivial ${\bf R}_{>0}^s$-bundle $\widetilde {\cal T}_{g}^s\to T_{g}^s$, where
the fiber over a point is the set of all $s$-tuples of horocycles in $F$, one horocycle about each puncture, and there is a 
corresponding trivial bundle $\widetilde{\cal T}_g^s/PMC=\widetilde {M}_{g}^s\to M_{g}^s$; the total spaces of
these bundles are the {\it decorated Teichm\"uller} and {\it decorated moduli spaces}, respectively, and they are studied in
[P1-P3].

\vskip .1in

\noindent In the bordered case $r\neq 0$, there are two geometric treatments of 
distinguished points in the boundary, which will be carefully described and compared in $\S$3 leading to two closely related
models of corresponding decorated moduli spaces.

\vskip .1in

\noindent 
The remainder of this section is dedicated to a discussion of the well-known canonical cell decomposition of decorated moduli space in
the punctured unbordered
case $r=0$ with $s> 0$ (as we tacitly assume in the remainder of this section),
which may be described either in the spirit of [P1-P3] (``in the hyperbolic setting'') or using quadratic differentials [St]
(``in the conformal setting'') as in [Ha] or [Ko].  Our discussion begins in the conformal setting, where moduli space is regarded as
the space of all equivalence classes of conformal structures under push-forward by any orientation-preserving diffeomorphism.

\vskip .1in

\noindent  A ``fatgraph'' or ``ribbon graph'' $G$ is a graph whose vertices have valence at least three plus the further structure of 
a cyclic ordering on the half-edges about each vertex.
A fatgraph $G$ may be ``fattened'' to a bordered surface as follows: start with disjoint neighborhoods in the plane of the vertices of
$G$, where the cyclic ordering is determined by the orientation of the plane;  glue orientation-preserving bands to these neighborhoods in
the natural way, one band for each edge of
$G$, to yield a topological surface $F_G\supseteq G$ with $G$ a spine of $F_G$.
The complement $F_G-G$ is a union of topological
annuli
$A_i\subseteq F_G$, $i=1,\ldots ,s\geq 1$, where each annulus $A_i$ has one boundary component $\partial _i'$ contained in $G$ and the
other
$\partial _i$ contained in
$\partial F_G$.

\vskip .1in

\noindent Letting $E(G)$ denote the set of edges of $G$, a
``metric'' on $G$ is any function $w\in ({\bf R}_{\geq 0} )^{E(G)}$.   
The curve $\partial _i'$ is a cycle
on $G$, and we define the ``length'' of $A_i$ to be $\ell _i (w)=\sum w(e)$, where the sum is over $e\in\partial _i'$ counted
with multiplicity.  A metric $w$ is said to be ``positive'' if $\ell _i (w)>0$ for each $i$, and we let $\sigma (G)$ denote the
space of all positive metrics on $G$.

\vskip .1in

\noindent Suppose that $w\in \sigma (G)$ is a positive metric with associated lengths $\ell _i>0$, for $i=1,\ldots ,s$.  We 
construct a metric surface homeomorphic to $F_g^s$ as follows, where $2-2g-s$ is the Euler characteristic of $G$.  Give each $A_i=A_i(w)$
the structure of a flat cylinder with circumference
$\ell _i$ and height unity, and isometrically identify these cylinders in the natural way along common
edges in
$\cup\{
\partial _i'\} _1^s$ as dictated by the metric and fatgraph. This yields a metric structure on $F_G$ so that the boundary component
$\partial _i$ of
$F_G$ is a standard circle of circumference $\ell _i$;
for each $i=1,\ldots ,s$, 
glue a standard flat disk of circumference $\ell _i$ to each $\partial _i$  and take
the center of the flat disk as a puncture $*_i$ to produce a
conformal structure on a surface which is identified with $F_g^s$. 

\vskip .1in

\noindent An analytic fact is:

\vskip .2in

\noindent {\bf Theorem A}~ [St;$\S$23.5]~~\it Suppose $w\in\sigma (G)$ is a positive metric on the fatgraph $G$, and let
$\ell _i=\ell _i (w)$ denote the corresponding lengths, for
$i=1,\ldots ,s$.  Then there is a unique meromorphic quadratic differential $q$ on $F_g^s$ so that for each $i=1,\ldots ,s$, 
the non-critical horizontal
trajectories of $q$ in $F$ foliate $A_i(w)\subseteq F_g^s$ by curves homotopic to the cores, and the residue of $\sqrt{q}$ at $*_i$ is
$\ell _i$.\rm

\leftskip=0ex\rm

\vskip .2in

\noindent Let $\mu _q$ denote the conformal structure on $F_g^s$ determined by the quadratic differential $q$.  Theorem~A thus describes
an effective construction $$(G,w )\mapsto \mu _q \times (\ell _i )_1^s\in{ T}_g^s\times {\bf R}_{>0}^s.$$

\vskip .1in

\noindent Fix a surface $F=F_g^s$, and consider the set $C_g^s$ of all homotopy classes of inclusions of fatgraph $G\subseteq
F$ where $G$ is a strong deformation retract of $F$. If $G\subseteq F$ and $w\in\sigma (G)$, then we may produce
another $G_w\subseteq F$ by contracting each edge $e\in E(G)$ with $w(e)=0$ to produce $G_w$. Identifying $E(G_w)$ with 
$\{ e\in E(G): w(e)\neq 0\}$ in the natural way, we thus induce
$w'\in \sigma (G_w)$ by requiring that $w'(e)=w(e)$ whenever $w(e)\neq 0$.

\vskip .1in

\noindent Consider 
$$U_g^s=\biggl [~\coprod _{(G\subseteq F)\in C_g^s}{\sigma (G)}~\biggr ]/\sim,$$
where $\coprod$ denotes disjoint union, and $(G^1,w ^1)\sim (G^2,w ^2)$ if and only if $G^1_{w_1}\subseteq F$ agrees with
$G^2_{w_2}\subseteq F$ as members of
$C_g^s$ and $w_1'\in\sigma (G^1)$ agrees with $w_2'\in\sigma (G^2)$. 
$PMC(F)$ acts on $U_g^s$ induced by $(\phi:G\to F)\mapsto (f\circ\phi :G\to F)$ for any diffeomorphism $f:F\to F$, and the quotient
``fatgraph complex''
is $G_g^s=U_g^s/PMC$; let $[G,w]\in G_g^s$ denote the class of $(G,w)\in U_g^s$.

\vskip .2in

\noindent{\bf Theorem B} [S;$\S$25.6] and [HM]~\it ~$(G,w )\mapsto \mu _q \times (\ell _i)_1^s$ induces real-analytic
homeomorphisms $
U_g^s\to { T}_g^s\times {\bf R}_{>0}^s$ and 
$G_g^s\to M_g^s\times {\bf R}_{>0}^s$.\rm

\vskip .2in

\noindent  By Theorem A, the mapping $G_g^s\to M_g^s\times {\bf R}_{>0}^s$ is well-defined and injective, and by Theorem B, this
mapping is surjective (but its inverse requires solving the Beltrami equation and is highly non-computable). 
There is thus a
$PMC$-invariant cell decomposition of $T_g^s\times{\bf R}_{>0}^s$ induced by the cell structure of $G_g^s$, which is the main result in
the conformal setting positing the existence of a cell decomposition of $M_g^s$.

\vskip .1in

\noindent In the hyperbolic setting, fix a (conjugacy class of) Fuchsian group $\Gamma$
uniformizing a point of $T_g^s$, and specify a collection of horocycles, one horocycle about each puncture
of
$F=F_g^s$ (called a ``decoration''); these data uniquely determine
a point $\tilde\Gamma\in \widetilde{\cal T}_g^s$ by definition.   The main ingredient for the cell decomposition of moduli space in the
hyperbolic setting is the  ``convex hull construction'' [P1;$\S$4], which assigns
to
$\tilde\Gamma\in \widetilde{\cal T}_g^s$ a corresponding point 
$(G_{\tilde\Gamma},w_{\tilde\Gamma})\in U_g^s$ as follows:

\vskip .1in

\noindent Via affine duality, we may identify the
open positive light-cone $L^+$ in Minkowski three-space with the collection of all horocycles in the hyperbolic plane 
(cf. $\S$2) and thus produce a
$\Gamma$-invariant set $B\subseteq L^+$ corresponding to the decoration, where we regard $\Gamma$ as acting via Minkowski
isometries.  In fact, $B$ is discrete in $L^+$, and we may consider the closed convex hull of $B$ in the vector space structure
underlying Minkowski space.  The extreme edges of the resulting
$\Gamma$-invariant convex body project to a collection of disjointly embedded arcs $\alpha _{\tilde\Gamma}$ connecting
punctures, and each component of $F-\cup{\alpha _{\tilde\Gamma}}$ is simply connected; we shall say that such an arc family
$\alpha_{\tilde\Gamma}$ ``fills'' $F$.  Given any arc family $\alpha$ filling $F$, we may define a subset
$$C(\alpha )=\{ \tilde\Gamma\in \widetilde{\cal T}_g^s:\alpha_{\tilde\Gamma}~{\rm is~homotopic~to}~\alpha \} .$$
The Poincar\'e dual of the cell decomposition $\cup\alpha _{\tilde\Gamma}$ of $F$ is a fatgraph $G$ embedded as a spine of $F$.
The explicit formula in terms of Minkowski geometry for the ``simplicial coordinate'' of an edge of $G$, gives a
positive metric $w_{\tilde\Gamma}$ on $G$; the simplicial coordinates are effectively the signed volume of the corresponding simplex
with light-like vertices, as we shall describe in $\S$2.

\vskip .1in

\noindent In the conformal setting, the effective construction maps $G_g^s\to (M_g^s\times{\bf R}_{>0}^s)$, while in the hyperbolic
setting, the effective construction (i.e., the convex hull construction) maps $\tilde M_g^s\to G_g^s$ in the opposite direction!
Thus, each of the hyperbolic and conformal treatments has its difficult theorem: surjectivity of the effective construction.
The proof of surjectivity in decorated
Teichm\"uller theory devolves to showing that the putative cells $C(\alpha )\subseteq\widetilde{\cal T}_g^s$ are in fact cells, and
this cellularity is proven in [P1;$\S$5] (independent of the Jenkins-Strebel theory) as described further in $\S$2. 

\vskip .1in

\noindent  
Just
as the conformal setting has a non-computable inverse $(M_g^s\times{\bf R}_{>0}^s)\to G_g^s$, there is a
non-computable (or at least, very difficult to compute) inverse $G_g^s\to \tilde M_g^s$ in the hyperbolic setting.  These ``arithmetic
problems'' are studied in [P2] and will be discussed in $\S$2; a solution to the arithmetic problem is thus a solution to an appropriate
Beltrami equation.

\vskip .1in

\noindent 
We believe that one may combine bounded distortion results of Epstein-Marden/Sullivan from the 1980's with arguments from [P1,$\S$6] 
and prove that the hyperbolic construction followed by the conformal construction has bounded distortion on open moduli space (and
the details of this will appear elsewhere).   For instance, one can thus view the arithmetic problem as the calculation of hyperbolic
geometry from conformal combinatorics, up to a bounded error.

\vskip .1in

\noindent It is worth emphasizing that the conformal and hyperbolic versions of the cell decomposition of $T_g^s$ agree combinatorially
but not pointwise.  Since $PMC$ acts more or less discretely on $T_g^s$, it is not surprising that a given cell decomposition can be
``jiggled'' to produce a combinatorially equivalent but pointwise distinct cell decomposition.  (On the other hand, another
hyperbolic version of the cell decomposition due to Bowditch-Epstein was shown in [P2] to agree {\sl pointwise} with the cell
decomposition based on the convex hull construction.)  Finally, 
there is further structure in the hyperbolic setting without analogue in the conformal setting as discussed in $\S$2-3 (from [P1] in the
punctured and [P4] in the bordered case).

\vskip .3in

\noindent {\bf 2. Coordinates on Teichm\"uller space for punctured surfaces}

\vskip .2in

\noindent There are several coordinatizations of the Teichm\"uller space ${\cal T}_g^s={\cal T}(F_g^s)$ of a punctured unbordered
surface which we shall next describe. All depend upon the Minkowski inner product $<\cdot ,\cdot >$ on ${\bf R}^3$ whose quadratic form
is given by
$x^2+y^2-z^2$ in the usual coordinates.   As is well-known, the upper sheet
${\bf H}=\{ u=(x,y,z)\in{\bf R}^3:<u,u>=-1~{\rm and}~z>0\}$ of the two-sheeted hyperboloid is isometric to the hyperbolic plane. 
Furthermore, the open positive light cone
$L^+=\{ u=(x,y,z)\in{\bf R}^3:<u,u> =0~{\rm and}~z>0\}$ is identified with the collection of all horocycles in ${\bf H}$ via 
the correspondence $u\mapsto h(u)=\{ w\in{\bf H}:<w,u>=-1\}$.

\vskip .1in

\noindent 
The first basic building block to be understood is a ``decorated geodesic'', by which we mean a pair of horocycles $h_0 ,h_1$ in the
hyperbolic plane with distinct centers, so there is a well-defined geodesic connecting the centers of $h_0$ and $h_1$; the two
horocycles may or may not be disjoint, and there is a well-defined signed hyperbolic distance $\delta$  between the horocycles
(taken to be positive if and only if 
$h_0\cap h_1=\emptyset$) as illustrated in the two cases of Figure~1.  The {\it lambda length} of the pair of horocycles $\{ h_0,h_1\}$ 
is defined to be the transform $\lambda (h_0,h_1)=\sqrt{2~{\rm exp}~\delta}$.  Taking this particular transform renders the
identification
$h$ geometrically natural in the sense that 
$\lambda (h(u_0),h(u_1))=\sqrt{-<u_0,u_1>}$, for $u_0,u_1\in L^+$ as one can check.

\vskip 1.8in

~~~{{{\epsffile{mosh1.epsf}}}}

\vskip .1in

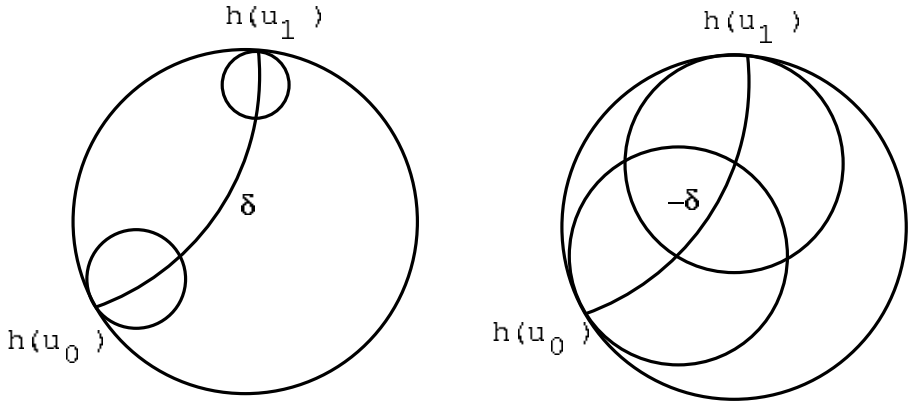
\centerline{{\bf Figure 1}~~{Decorated geodesics.}}

\vskip .1in

\noindent By an {\it arc family} in $F_g^s$, we mean the isotopy class of a family of essential arcs disjointly embedded in $F_g^s$
connecting punctures, where no two arcs in a family may be homotopic rel punctures.  If $\alpha$ is an arc family so that each component
of
$F-\cup \alpha$ is a polygon, then we say that $\alpha$ {\it fills} $F$, and in the extreme case that $\alpha$ is maximal so that each
complementary component is a triangle, then $\alpha$ is called an {\it ideal triangulation}.  

\vskip .1in
 \noindent Let us choose an ideal triangulation $\tau$
of $F_g^s$ once and for all and regard this as a ``choice of basis'' for the first global coordinatization:

\vskip .2in

\noindent {\bf Theorem 1} [P1; Theorem 3.1]\it ~For any ideal triangulation $\tau$ of $F$, the assignment of lambda lengths defines a
real-analytic homeomorphism $$\widetilde {\cal T}_g^s\to{\bf R}_{>0}^{\tau}.$$\rm

\vskip .2in

\noindent  For the proof, we refer the reader to [P1], which gives an effective construction of the
underlying Fuchsian group from the lambda lengths on $\tau$ (and is the decorated version of an easy case of Poincar\'e's original
construction of Fuchsian groups).

\vskip .1in

\noindent The next basic building block to be understood is a ``decorated triangle'', by which we mean a triple $h_0,h_1,h_2$ of
horocycles with distinct centers, so these centers determine an ideal triangle in the hyperbolic plane.  We shall refer to a
neighborhood of a cusp of this triangle as a ``sector'' of the triangle, so there is one horocylic segment contained in each sector, and
we define the {\it h-length} of the sector to be half the hyperbolic length of the corresponding horocyclic segment; see Figure~2, where
we have identified an arc with its lambda length for convenience.  One can think of the h-length as a kind of angle
measurement between geodesics asymptotic to a puncture.

\vskip 1.8in

\hskip .6in{{{\epsffile{mosh2.epsf}}}}

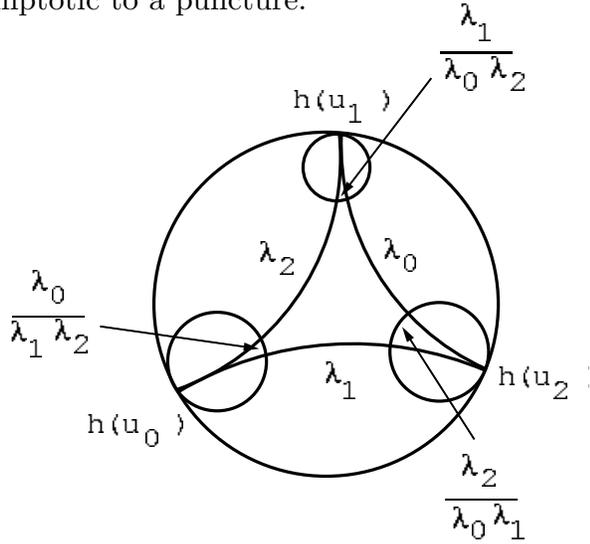
\centerline{{\bf Figure 2}~~{Decorated triangle.}}

\vskip .1in

\noindent {\bf Lemma 2}~~~\it Fix a decorated ideal triangle $h_0,h_1,h_2$, and define the lambda lengths
$\lambda _i=\lambda
(h_j,h_k)$, where
$\{ i,j,k\}=\{ 0,1,2\}$.  Then 

\vskip .1in

\noindent {\bf a)} {\rm [P1;Proposition 2.8]} As in Figure~2, the h-length of the sector opposite $\lambda _i$ is given by
${{\lambda _i}\over{\lambda _j\lambda _k}}~~{\rm
for}~~\{ i,j,k\}=\{ 0,1,2\}$.

\vskip .1in

\noindent {\bf b)}  {\rm [P2;Propostion 2.3]} There is a point in the hyperbolic plane equidistant to $h_0,h_1,h_2$ if and only if
$\lambda _0,\lambda _1,\lambda _2$ satisfy all three possible strict triangle inequalities, and this equidistant point is unique in
this case as illustrated in Figure~3; furthermore, for a fixed ideal triangle, all points in the hyperbolic plane arise for some
decoration on it.

\vskip .1in

\noindent {\bf c)} {\rm [P5;WP volumes;Theorem~3.3.6]} The Weil-Petersson K\"ahler two-form is given by
$$-2\sum~d{\rm log}~\lambda _0\wedge d{\rm log}~\lambda _1+
d{\rm log}~\lambda _1\wedge d{\rm log}~\lambda _2+
d{\rm log}~\lambda _2\wedge d{\rm log}~\lambda _0,$$
where the sum is over all triangles complementary to $\tau$ with edges $\lambda _0,\lambda _1,\lambda _2$ in this
cyclic order in accordance with the orientation of $F$.

\rm

\vskip .1in

\noindent The proofs of a) and b) follow from direct calculation, and c) follows from Wolpert's formula plus calculation.

~~\vskip 1.8in

\hskip .1in{{{\epsffile{mosh3.epsf}}}}

\hskip .45in{inside triangle}\hskip 1.5in{outside triangle}

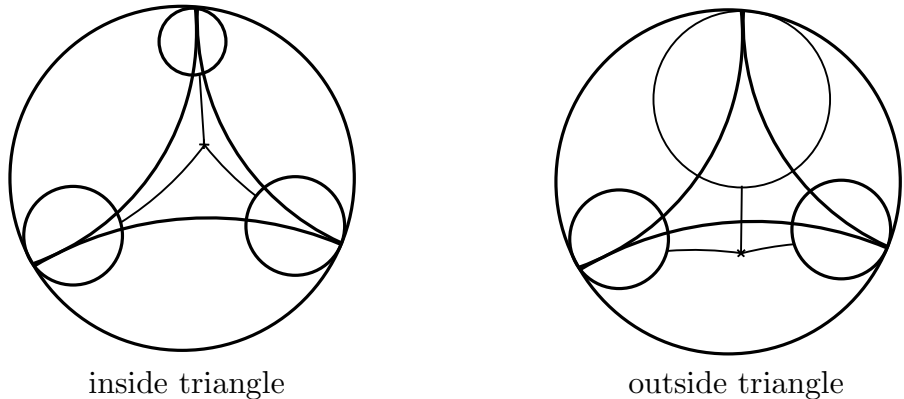
\centerline{{\bf Figure 3}~~{Equidistant points to triples of horocycles.}}

\vskip .2in

\noindent The last basic building block to be understood is a ``decorated quadrilateral'' by which we mean four horocycles with
distinct centers, which thus determine an ideal quadrilateral.  Let us furthermore triangulate this quadrilateral by adding one of
its diagonals $e$ and adopt the notation of Figure~4 for the nearby lambda lengths, where again we identify an arc with its lambda
length for convenience.  Furthermore, adopt the notation of Figure~4 for the nearby h-lengths.  Define the {\it simplicial coordinate}
of the edge $e$ to be
$$E={{a^2+b^2-e^2}\over{abe}}+{{c^2+d^2-e^2}\over{cde}}=(\alpha +\beta-\epsilon)+(\gamma +\delta -\phi ).$$
In fact, $abcdE$ is the signed volume of the corresponding tetrahedron in Minkowski space, as one can check.

\vskip 1.8in

\hskip .9in{{{\epsffile{mosh4.epsf}}}}

\vskip .2in

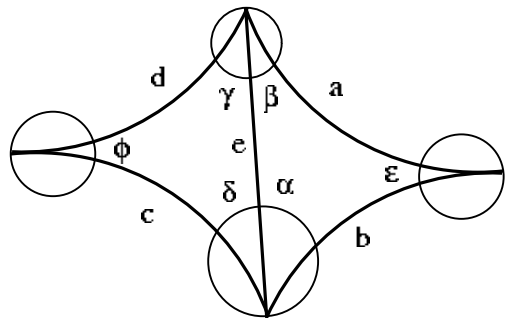
\centerline{{\bf Figure 4}~~{Decorated quadrilateral.}}

\vskip .2in

\noindent {\bf Lemma 3} [P1:$\S$2]~\it Adopting the notation of Figure~4, we have

\vskip .1in

\noindent {\bf a)}~the {\rm coupling equation} holds: $\alpha\beta =\gamma\delta$;

\vskip .1in

\noindent {\bf b)}~the cross-ratio of the four ideal vertices is given by ${{ac}\over{bd}}$;

\vskip .1in

\noindent {\bf c)}~the diagonal $f$ of the quadrilateral other than diagonal $e$ has lambda length given by a {\rm Ptolemy
transformation}:
$f={{ac+bd}\over e}$.\rm

\vskip .2in

\noindent The proofs follow from easy calculation.  The coupling equations give a variety determined by a 
collection of coupled quadric equations in the vector space spanned by the
sectors, where the simpicial coordinates are given by linear constraints, and this variety is identified with decorated Teichm\"uller
space.

\vskip .1in

\noindent Fix an ideal triangulation $\tau$ of $F_g^s$.  Define a {\it cycle of triangles} $(t_j)_1^n$ to be a collection of 
triangles complementary to $\tau$ so that
$t_j\cap t_{j+1}=e_j$, for each $j=1,\ldots ,n$, and the index $j$ is
cyclic (so that
$t_{n+1}=t_1$).  If the edges of $t_j$ are $\{ e_j,e_{j-1},b_j\}$, for $j=1,\ldots ,n$, then the collection $\{
b_j\} _1^n$ is called the {\it boundary} of the cycle, and the edges $\{ e_j\}$ are called the 
{\it consecutive edges} of the cycle.  A function $X:\{ {\rm edges~of}~\tau\} \to{\bf R}_{\geq 0}$ satisfies the 
{\it no vanishing cycle} condition if for any
cycle of triangles with consecutive edges $e_{n+1}=e_1,\ldots ,e_n$, we have $\sum _{j=1}^n X(e_j)>0$.

\vskip .1in

\noindent A nifty little fact is the ``telescoping property'' of simplicial coordinates, as follows:  If $e_1,\ldots ,e_n$ are the
consecutive edges of a cycle of triangles with boundary $b_1,\ldots ,b_n$, and
the sector in
$t_i$ opposite
$b_i$ has h-length
$h_i$, then
$\sum _{i=1}^n E_i=2\sum _{i=1}^n h_i$, where
the simplicial coordinate of edge $e_i$ is given by $E_i$.  The proof follows immediately from the definition of simplicial
coordinates.  

\vskip .1in

\noindent 
Suppose that $\tilde \Gamma\in\widetilde{\cal T}_g^s$ gives rise via the convex hull construction to
the ideal cell decomposition $\alpha _{\tilde\Gamma}$.  Let us complete $\alpha _{\tilde\Gamma}$ to an ideal triangulation $\tau$
in any manner (triangulating the polygonal complementary regions).
From the very definition of the convex hull construction and simplicial coordinate, we find that the simplicial coordinate of each
edge is non-negative (by convexity) and satisfies the no vanishing cycle condition (since complementary regions are polygonal).

\vskip .1in

\noindent Conversely, the main hard result of decorated Teichm\"uller theory is

\vskip .2in

\noindent {\bf Theorem 4}~  [P1;Theorem 5.4]\it Given an ideal triangulation $\tau$ and any $X_e\in{\bf R}_{\geq 0}$, for each
$e\in\tau$,  with no vanishing cycles, there is a unique corresponding assignment of lambda lengths in ${\bf R}_{>0}^\tau$ determining
a decorated hyperbolic structure
$\tilde\Gamma\in\widetilde{\cal T}_g^s$ so that $\alpha _{\tilde\Gamma}=\tau -\{ e\in\tau : X_e=0\}$, and the lambda lengths
on $\tau$ realize the putative simplicial coordinates $X_e$.\rm

\vskip .2in

\noindent The idea of the proof is to consider the vector space
spanned by the sectors of an ideal triangulation $\tau$ and introduce the ``energy'' 
$$K=\sum  \biggl [ {\rm log }~ \bigl ({{\alpha\beta}\over{\gamma\delta}}\bigr )\biggr ]^2,$$
where the sum is over all edges $e$ of $\tau$ in the notation of Figure~4.
In fact, solving for lambda lengths amounts to minimizing $K$ subject to the constraints given by the simplicial coordinates, the
constrained gradient flow of $K$ is dissipative with a unique minimum, this unique constrained minimum of $K$ solves the coupling
equations, and this gives the desired lambda lengths from $\alpha \beta =e^{-2}=\gamma\delta$.

\vskip .1in

\noindent Thus, the explicit formulas for simplicial coordinates in terms of lambda lengths can be uniquely inverted provided the
simplicial coordinates are non-negative and satisfy the no vanishing cycle condition.  This explicit calculation is called the 

\vskip .2in

\noindent {\bf Arithmetic Problem}~[P2]~\it Fix any ideal triangulation of any punctured surface, and calculate the lambda lengths from
the simplicial coordinates.\rm

\vskip .2in

\noindent {\bf Lemma 5}~[P1;Lemma 5.2]~\it The non-vanishing cycle condition on simplicial coordinates implies that the lambda lengths
on any complementary triangle satisfy all three strict triangle inequalities; in particular, the product of the lambda length and the
simplicial coordinate of an edge is bounded above by four.\rm

\vskip .2in

\noindent {\bf Proof}~ 
Adopt the notation in the definition of simplicial coordinates for the lambda lengths near an edge $e$.  If
$c+d\leq e$, then $c^2+d^2-e^2\leq -2cd$, so the non-negativity of the simplicial coordinate $E$ gives
$0\leq cd  [(a-b)^2-e^2]$, and we find a second edge-triangle pair so that the triangle inequality fails.
This is a basic algebraic fact about simplicial coordinates.
It follows that
if there is any such
triangle $t$ so that the triangle inequalities do fail for $t$, then there must be a cycle of triangles of such failures. 
This possibility is untenable since if
$e_{j+1}\geq b_j+e_j$, for $j=1,\ldots n$, where we again identify an arc with its lambda length, then upon summing and canceling
like terms, we find $0\geq \sum _{j=1}^n b_j$, which is absurd since lambda lengths are positive.

\vskip .1in

\noindent The second part follows from the elementary observation in the notation of Figure~4 that
$$eE={{a^2+b^2-e^2}\over{ab}}+{{c^2+d^2-e^2}\over{cd}},$$
i.e., each summand is twice a cosine, and hence $eE\leq 4$.~~\hfill{\it q.e.d.}

\vskip .2in

\noindent By taking Poincar\'e duals of filling arc families, there is established an isomorphism between fatgraphs and
ideal cell decompositions, and the simplicial 
coordinates are the analogues of the metric lengths of edges in the conformal treatment.  Indeed, we may restate the main theorem of
$\S$1 in the hyperbolic setting using arc families as:

\vskip .2in

\noindent {\bf Theorem 6}~[P1;Theorem 5.5]\it  ~There is a $PMC(F_g^s)$-invariant cell decomposition of $\widetilde{\cal T}_g^s$
isomorphic to the geometric realization of the partially ordered set of all arc families filling $F_g^s$.\rm

\vskip .4in

\noindent {\bf 3. Bordered surfaces}\vskip .2in

\noindent We next turn our attention to bordered surfaces.  Not only
does much of the decorated Teichm\"uller theory extend naturally to this setting, but there are also new aspects of the theory to be
described.

\vskip .1in

\noindent  
In this section, we shall assume only that
$r\neq 0$ (and allow both cases $s=0$ or $s>0$).  Enumerate the (smooth) boundary components of $F$ as $\partial _i$, where $i=1,\ldots
,r$, and set
$\partial=\cup\{
\partial _i\} _1^r$.  
\vskip .1in

\noindent Let $Hyp(F)$ be the space of all hyperbolic metrics with geodesic boundary on
$F$.  Define the {\it moduli space} to be
$$M=M(F)= \bigl [Hyp(F)\times (\prod _{1}^r\partial _i)\bigr ]/\sim ,$$
where $\sim$ is the equivalence relation generated by push-forward of metric under orientation-preserving diffeomorphism
$$f_*\bigl (\Gamma ,(\xi _i )_1^r)\bigr ) = (f_*(\Gamma ),(f(\xi _i)_1^r),$$ where $\Gamma\mapsto f_*(\Gamma )$ is the usual
push-forward of metric on $Hyp(F)$.

\vskip .2in

\noindent Let $\ell _i(\Gamma)$ denote the hyperbolic length of $\partial _i$ for $\Gamma\in Hyp(F)$, and define the first model 
$\tilde M=\tilde M(F)$ of {\it
decorated moduli space} to be
$$\bigl\{ \bigl (\Gamma ,(\xi _i )_1^r,(t_i)_1^r\bigr ):\Gamma\in Hyp(F), \xi _i\in\partial _i, 0<t_i< \ell
_i(\Gamma ), i=1,\ldots ,r\bigl\} /\sim ,$$
where $\sim$ denotes push-forward by diffeomorphisms on $(\Gamma ,(\xi _i )_1^r )$, as before, extended by the trivial action on $(t_i
)_1^r$.  A point of $\tilde M$ is thus represented by $\Gamma\in Hyp(F)$ together with a pair of points $\xi _i\neq p_i$ in each
component $\partial _i$, where $p_i$ is the point at hyperbolic distance $t_i$ along $\partial _i$ from $\xi _i$ in the natural
orientation. (In the special case
$g=0=s=r-2$ where
$F$ is an annulus, we define
$\tilde M(F)$ to be the collection of all configurations of two distinct labeled points in an abstract circle of some radius.)

\vskip .1in

\noindent Turning to the second model for decorated moduli space, begin with a smooth surface $F$ with smooth boundary, choose one
distinguished point
$d_i\in\partial _i$ in each boundary component, and set
$D=\{ d_i\} _1^r$.  
Define a {\it quasi hyperbolic
metric} on
$F$ to be a hyperbolic metric on
$F^\times =F-D$ so that $\partial _i^\times =\partial _i-\{ d_i\}$ is totally geodesic, for each $i=1,\ldots ,r$, and set $\partial
^\times=\cup\{
\partial _i^\times\} _1^r$.  
To illuminate the definition, take a hyperbolic metric on a  once-punctured annulus $A$ with geodesic boundary and the simple
essential geodesic arc
$a$ in $A$ which is asymptotic in both directions to the puncture; the induced metric on a
component of $A-a$ gives a model for the quasi hyperbolic metric near a component of $\partial ^\times $ as in Figure~5.  

~~\vskip 2.8in

\hskip 1.9in{\epsffile{sphfig1.epsf}}

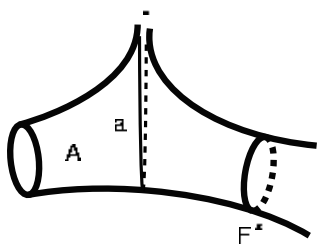
\centerline{{\bf Figure 5}~~{The model for $F^\times$.}}

\vskip .2in

\noindent In this second model, there is a 
version of {\it decorated Teichm\"uller space} as the space
$\widetilde{\cal T}=\widetilde{\cal T}(F)$ of all quasi hyperbolic metrics on $F-D$, 
where we furthermore specify for each $d_i$ a horocyclic segment
centered at $d_i$ (again called a ``decoration''), modulo push-forward by diffeomorphisms of $F-D$ which are isotopic to the identity,
where the action is trivial on hyperbolic lengths of horocylic segments.  In fact, $\widetilde{\cal T}$ is homeomorphic to an open ball of
dimension
$6g-6+5r+2s$. The second model for the {\it decorated moduli
space} is
$\widetilde{\cal M}=\widetilde{\cal M}(F)=\widetilde{\cal T}(F)/PMC(F)$.

\vskip .1in

\noindent For a fixed quasi hyperbolic metric on $F$, there is a unique separating geodesic $\partial _i^*\subseteq F$ in the homotopy
class of $\partial _i$ .  Excising from $F-\cup\{\partial _i^*\}_1^r$ the components containing
points of $\partial ^\times$ yields a surface $F^*$, which inherits a hyperbolic metric with geodesic boundary.  (In the special
case that $F$ is an annulus, the surface $F^*$ collapses to a circle.)  $\tilde\Gamma\in\widetilde{\cal T}$
has its underlying hyperbolic metric given by a conjugacy class of Fuchsian group $\Gamma$ for
$F^*$.

\vskip .1in

\noindent Define an {\it arc} in $F$ to be a smooth path $a$
embedded in $F$ whose endpoints lie in $D$ and which meets $\partial F$ transversely, where
$a$ is not isotopic rel its endpoints to a path lying in $\partial F -D$.  Two arcs
are {\it parallel} if there is an isotopy between them
fixing $D$ pointwise.  
An {\it arc family} in
$F$ is the isotopy class of a set of disjointly embedded
arcs in
$F$, no two components of which are parallel.  
A collection $\alpha$ of arcs representing an arc family in $F$ is said to {\it meet efficiently} an 
arc or curve
$C$ embedded in
$F$ if there are no bigons in $F$ complementary to $\alpha\cup C$.

\vskip .1in

\noindent If $\alpha$ is an arc family in $F$ where each component of
$F-\cup \alpha $ is a polygon or a once-punctured polygon, then we say that
$\alpha $ {\it quasi fills} $F$.  In the extreme case that each component is either a triangle
or a once-punctured monogon, then
$\alpha$ is called a {\it quasi triangulation}.

\vskip .2in

\noindent{\bf Theorem 7} ~[P4;Theorem 1]~~\it For any quasi triangulation $\tau$ of $F$, the assignment of lambda lengths
defines a real-analytic homeomorphism
$$\widetilde{\cal T}_{g,r}^s\to {\bf R}_{>0}^{\tau\cup\partial ^\times}.$$\rm

\vskip .2in

\noindent The proof is verbatim that of [P1; Theorem~3.1] briefly mentioned before.  In fact, much of the number-theoretic utility of
lambda lengths in the punctured case [P1-P3] extends without great difficulty to the bordered case as discussed in further detail in the
Side-Remark after Theorem 1 of [P4].

\vskip .1in

\noindent One introduces simplicial coordinates in the bordered case in analogy to the previous discussion, where if
$e\in\partial ^\times$, then it bounds a triangle on only one side, say with edges $a,b,e$, and we define its simplicial coordinate to
be
$E=2~{{a^2+b^2-e^2}\over{abe}}$.  

\vskip .1in

\noindent Fix a quasi triangulation $\tau$ of $F$, and define the subspace
$$\eqalign{
\tilde C(\tau )&=\{ (\vec{y},\vec{x})\in{\bf R}^{
\partial ^\times}\times ({\bf R}_{\geq 0} )^{\tau }:
~{\rm there~are~no~vanishing~cycles~or~arcs}\}\cr
&\subseteq {\bf R}^{\tau\cup\partial ^\times} ,\cr
}$$
where the coordinate functions are the simplicial coordinates (rather than the lambda lengths as in the previous theorem).
We say that there are ``no vanishing cycles'' as before for a point of $C(\tau )$ if there is no essential simple closed
curve $C\subseteq F$ meeting a representative $\tau$ efficiently so that
$$0=\sum _{p\in C\cap\cup\tau} x_p,$$ where $p\in C\cap a$, for
$a\in\tau $, contributes to the sum the coordinate $x_p$ of
$a$.  We say that there are ``no vanishing arcs'' for a point of $C(\tau )$ if there is no essential simple arc
$A\subseteq F$ meeting $\tau$ efficiently and properly embedded in $F$
with endpoints disjoint from $D$
so that
$$0=\sum _{p\in A\cap \partial ^\times} y_p~~+~~\sum
_{p\in A\cap\tau} x_p,$$
where again $x_p$ and $y_p$ denote the coordinate on $a$
at an intersection point $p=A\cap a$ or $p=C\cap a$ for $a\in\tau$.  
This is a convex constraint on $\vec{y}$ for fixed 
$\vec{x}$.

\vskip .1in

\noindent Just as in the case of punctured unbordered surfaces, $\tilde\Gamma$ gives rise to a quasi filling arc family $\alpha
_{\tilde\Gamma}$ via the convex hull construction; fixing the topological type of $\alpha
_{\tilde\Gamma}$ and varying
$\tilde{\Gamma}$ then gives a cell in a natural decomposition of $\widetilde{\cal T}$. 

\vskip .2in

\noindent {\bf Theorem 8} ~[P4; Theorem 3]~~\it A point $\tilde\Gamma\in\widetilde{\cal T}$ gives rise to a quasi filling arc family
$\alpha _{\tilde\Gamma}$ via the convex hull construction as well as a tuple of simplicial coordinates $(\vec y,\vec x)\in\tilde
C(\tau)$ for any quasi triangulation
$\tau\supseteq\alpha_{\tilde\Gamma}$, where $\vec x$ vanishes on $\tau-\alpha _{\tilde\Gamma}$.
Furthermore, this gives a real-analytic homeomorphism of the
decorated Teichm\"uller space of
$F$ with $\bigl [\bigcup _\tau \tilde C(\tau )\bigr ]/\sim$, where 
$(\tau _1,\vec y_1,\vec x_1)\sim (\tau _2,\vec y_2,\vec x_2)$ if  
$\vec x_1$ agrees with $\vec x_2$ on 
$\tau _1-\{a\in\tau _1:x^1=0\}=\tau _2-\{a\in\tau _2:x^2=0\}$,
where $x^j$ denotes the $\vec x$ coordinate on $a$, for $j=1,2$, and $\vec y_1=\vec y_2$.
\rm

\vskip .2in

\noindent 
Define a
variant of the decorated Teichm\"uller space
$$\widehat{\cal T}=\widehat{\cal T}(F)=\bigl  [\bigcup _\tau \hat C(\tau )\bigr ]/\sim,$$
where $\sim$ is as in Theorem~8, and for membership in $\hat C(\tau )$ we demand that there are no vanishing cycles
or arcs and for any triangle
$t\subseteq F$ complementary to
$\tau$, and the lambda lengths on the edges of $t$ satisfy all three possible strict triangle inequalities.  Using cycles of
triangles, the quotient
$\widehat{\cal M}=\widehat{\cal T}/PMC$ can be shown [P4;Lemma 5] to be a strong deformation retract of $\widetilde{\cal
M}$.  Using the existence of equidistant points to suitable triples of horocycles as in Lemma~2b in order to determine distinguished
points, one can then prove

\vskip .2in

\noindent{\bf Theorem 9}~[P4; Theorem 8] \it There is a real-analytic homeomorphism between $\widehat{\cal M}$ and $\tilde M$.\rm

\vskip .2in

\noindent Thus, our two models $\tilde M$ and $\widetilde{\cal M}$ of decorated moduli space are homotopy equivalent (and indeed
are homeomorphic except for the required passage to $\widehat{\cal M}\subseteq \widetilde{\cal M}$).

\vskip .4in

\noindent{\bf 4. The arc complex of a bordered surface}

\vskip .2in

\noindent Let us inductively build a simplicial complex $Arc' (F)$, where
there is one
$p$-simplex
$\sigma(\alpha )$ for each arc family $\alpha $ in
$F$ of cardinality $p+1$; the simplicial structure of $\sigma (\alpha )$ is the natural one, where
faces of 
$\sigma (\alpha )$ correspond to sub arc families of $\alpha $.  Identifying the open standard $p$-simplex with the
collection of all real-projective $(p+1)$-tuples of positive reals assigned to the vertices, $Arc'(F)$ is identified with the
collection of all arc families in
$F$ together with a real-projective weighting of non-negative real numbers, one such number assigned to each component of
$\alpha$.

\vskip .1in

\noindent $PMC(F)$ acts on $Arc'(F)$ in the natural way,
and we define the {\it arc complex} of $F$ to be the quotient $$Arc(F)=Arc'(F)/PMC(F).$$
If $\alpha $ is an arc
family in $F$ with corresponding simplex $\sigma (\alpha )$ in $Arc'(F)$, then we shall let $ [\alpha
]$ denote the
$PMC(F)$-orbit of $\alpha $ and $\sigma  [\alpha ]$ denote the quotient of $\sigma (\alpha )$
in $Arc(F)$.  A subspace of $Arc(F)$ of special interest is the space of quasi-filling arc families with positive projective weight
$$Arc_\#(F)=\{ \sigma [\alpha ]: \alpha ~{\rm quasi-fills}~F\} .$$

\vskip .1in

\noindent Let $M(F)/{\bf R}_{>0}$ denote the quotient of the moduli space of the bordered surface $F$ by the ${\bf R}_{>0}$-action
by homothety on the tuple of hyperbolic lengths of the geodesic boundary components.  A further analysis of the homeomorphism in
Theorem~9 shows

\vskip .2in

\noindent {\bf Theorem 10} [P4;Theorem 14] \it For any bordered surface $F\neq F_{0,2}^0$, $Arc_\# (F)$
is proper homotopy equivalent to
$M(F)/{\bf R}_{>0}$.\rm

\vskip .3in

\noindent  Furthermore, each of $M(F)/{\bf R}_{>0}$ and $Arc_\# (F)$ admit natural $(S^1)^r$-actions (moving distinguished points in the
boundaries for the former and ``twisting'' arc families around the boundary components for the latter); the proper homotopy equivalence
in Theorem 10 is in fact a map of $(S^1)^r$-spaces.  This is the main result of [P4].  It follows from Theorem 10 plus the sphericity
conjecture 12 discussed in the next section that $M(F)/{\bf R}_{>0}$ compactifies to a sphere.

\vskip .2in

\noindent Notice that the space of all projectivized laminations ${\cal PL}(F)$ is itself a sphere (cf. [PH]) which in particular
contains the complex $Arc'(F)$.  It is well-known that ${\cal PL}(F)/PMC(F)$ is a space whose largest Hausdorff quotient is a
singleton, yet it contains as a dense open set the quotient $Arc(F)$, whose further study we undertake in the next section.

\vskip .4in

\noindent {\bf 5. Sphericity}\vskip .2in

\noindent Extending the notation from before, let $F=F_{g,\vec\delta }^s$ denote a fixed smooth
surface with genus $g\geq 0$, $s\geq 0$ punctures, and $r\geq 1$ boundary components, where the $i^{\rm th}$ boundary component
comes equipped with $\delta _i\geq 1$ distinguished points $d_i$ as well, for $i=1,\ldots ,\delta _i$ and $\vec\delta
=(\delta _1,\delta _2,\dots
\delta _r)$.

\vskip .1in

\noindent Let $Arc(F)=Arc '(F)/PMC(F)$ denote the arc complex of $F$ as defined before, and recall the following combinatorial fact.

\vskip .2in

\noindent {\bf Classical Fact 11} \it For $F=F_{0,(n)}^0$, $Arc(F)=Arc_\# (F)=Arc'(F)$ is PL-homeomorphic to the sphere of dimension
$n-4$ provided $n\geq 4$.\rm

\vskip .2in

\noindent Gian-Carlo Rota told us that sphericity for polygons was almost certainly known to Whitney.  This is a special
case of the following general conjecture. 

\vskip .2in

\noindent {\bf Sphericity Conjecture 12}~[P3;Conj.B]\it ~~The arc complex $Arc(F)$ is PL-homeomorphic to the sphere of dimension
$6g-7+3r+2s+\delta _1+\delta _2+\cdots +\delta _r$ provided this number is non-negative.\rm

\vskip .1in

\noindent This result has been a long-standing goal for us.  We had proved this conjecture for $g=0=r-1$ in [P3] and shall recall this
proof and add all the details here.  New tools have arisen to address the general sphericity conjecture in collaboration with Dennis
Sullivan.

\vskip .1in

\noindent   Let
$$D=D(F)=\cup\{
d_i:i=1,2,\dots ,r\}$$ denote the set of all the distinguished points in the boundary
of $F$, 
$$\Delta =\Delta (F)=\delta _1+\delta _2+\cdots +\delta _r,$$ be the cardinality of $D$, and set 
$$N=N(F)=6g-7+3r+2s+\Delta$$
to be the dimension of the arc complex.
Furthermore, in the special case that $\Delta =r$ so $\vec\delta$ is an $r$-dimensional vector each of whose
entries is unity, then we shall also sometimes write $F_{g,r}^s=F_{g,\vec\delta}^s$.  Notice that $N\geq 0$ except for the
triangle $F_{0,(3)}^0$ and once-punctured monogon $F_{0,1}^1=F_{0,(1)}^1$.

\vskip .1in

\noindent It is convenient in the sequel to imagine weighted arcs in the spirit of train tracks [PH] as follows.
If $\alpha $ is an arc family, then a {\it weighting} on $\alpha $ is the assignment of 
a non-negative real number, the {\it weight}, to each arc comprising $\alpha $, and the 
weighting is {\it positive} if each weight is positive.  A {\it projective weighting} is the
projective class of a weighting on an arc family.  Since the points of a simplex are described by 
``projective non-negative weightings on its vertices'', we
may identify a point of
$\sigma (\alpha )$ with a projective weighting on $\alpha $.  The
projective positive weightings correspond to the interior $Int~{\sigma }(\alpha )$ of $\sigma 
(\alpha )$, and
$Arc'(F)$ itself is thus identified with the collection of all projective positively weighted arc families
in $F$.  

\vskip .1in
 \noindent Suppose that $\alpha =\{ a_0,\ldots ,\alpha _p\}$ is an arc family.  A useful realization of a positive weighting
$w=(w_i)_0^p$ on the respective components of $\alpha$ is as a collection of $p+1$ ``bands'' $\beta _i$
disjointly embedded in the interior of $F$ with the width of $\beta _i$ given by $w_i$, where $a_i$ traverses the length of $\beta
_i$, for $i=0,\ldots ,p$.  The bands coalesce near the boundary of $F$ and are attached one-to-the-next in the natural manner
to form ``one big band'' for each of the $\Delta$ distinguished points as illustrated in Figure~6.  (More precisely, we are
constructing here a ``partial measured foliation from a weighted arc family regarded as a train track with stops''; for further
details, see [PH] or [KLP;$\S$1].)

~~\vskip 1.8in

\hskip .3in{{{\epsffile{band.epsf}}}}

\vskip -.2in

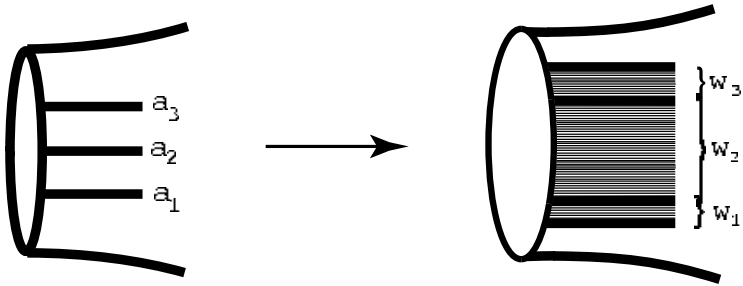
\centerline{{\bf Figure 6}~~{The band picture of a weighted arc family.}}

\vskip .1in

\noindent Given an arc family $\alpha $ in $F$, consider the surface $F-\cup\alpha$.  The non-smooth points or 
{\it cusps} on the frontiers of the components of $F-\cup\alpha $ give rise to distinguished points in the
boundary of
$F-\cup\alpha $ in the natural way, and this surface together with these distinguished points is
denoted
$F_{\alpha }$.  Our basic approach to sphericity is by induction on $N(F)>N(F_\alpha )$, and more will subsequently be said
on this point.

\vskip .1in

\noindent We begin the discussion of sphericity with a number of examples.

\vskip .2in\leftskip .2in

\noindent {\bf Example 1}~ If $g=s=0$ and $r=1$, then $F$ is a polygon with
$\delta _1$ vertices, $N=\delta _1-4$, $PMC(F)$ is trivial, and we are in the setting of the Classical Fact.  When
$\delta _1=4$ with
$F$ a quadrilateral, we have $N=0$, and there are exactly two chords, which cannot be simultaneously disjointly
embedded.  Thus,
$Arc(F)=Arc'(F)$ consists of two vertices, namely, $Arc(F)$ is the 0-sphere.  (The reader may likewise directly
investigate the case
$N=1$ where
$F$ is a pentagon, and the usual ``pentagon relation'' shows that $Arc(F)$ itself is also a pentagon.)  

\leftskip=0ex

\vskip .2in

\noindent More generally, in order to
handle all polygons inductively, we recall

\vskip .2in

\noindent {\bf Lemma 13}~~\it If $F'$ arises from $F$ by specifying one extra distinguished point in the
boundary of
$F$, then
$Arc(F')$ is PL-homeomorphic to the suspension of $Arc(F)$ provided $Arc(F)\neq\emptyset$.\rm

\vskip .1in

\noindent {\bf Proof}~ If $p$ is the
extra distinguished point of $F'$ in the boundary of $F$, let $p_0,p_1,p,p_2$ 
be the consecutive distinguished
points lying in a single boundary component of $F'$ in its induced orientation, where the $p_j$
may not be distinct.   There is a natural
mapping of 
$Arc'(F)\times [0,1]$ to
$Arc'(F')$ using partial measured foliations gotten by sliding proportion $t\in [0,1]$ of the band that hits $p_1$ over to
$p$, which descends to an embedding of
$(Arc(F)\times [0,1])/\sim)$ to $Arc(F')$; the equivalence relation
$\sim$ collapses the interval $[0,1]$ to a point for any arc family that 
does not hit $p_1$. There is an arc 
$a_1$ near the boundary which connects
$p_1$ to $p_2$ and another arc $a_0$ near the boundary which connects 
$p_0$ and $p$, each of which is essential in $F'$
but either non-essential or non-existent in $F$.  

\vskip .1in

\noindent The surface $F'_{\{ a_i\}}$ has two components, one of which is a triangle,
and the other of which is homeomorphic to $F$.  The 
collection $C_i$ of projectively weighted arc families in
$F'$ which either contain or to which we may add $a_i$ is thus homeomorphic to the topological join
$C_i\approx [a_i]*Arc(F)$ of $Arc(F)$ with a point, namely, the vertex $\sigma [\{ a_i\} ]$ corresponding to $a_i$, for
$i=0,1$.

\vskip .1in

\noindent Identifying the copy of $Arc(F)$ lying in $C_i$ with $$\bigl [(Arc(F)\times \{ i\} )/\sim \bigr ]~~\subseteq
~~\bigl [(Arc(F)\times
\{ 0,1\} )/\sim\bigr ]$$ in the natural way for $i=0,1$, 
we find a mapping of the 
suspension of $Arc(F)$ onto $Arc(F')$ which is
evidently a PL-homeomorphism.\hfill{\it q.e.d.}

\vskip .2in\leftskip=0ex

\noindent {\bf Proof of the Classical Fact} Induct on the number $\delta _1\geq 4$
of distinguished points using Lemma~13 starting from the basis step $\delta _1=4$  treated
in Example~1 (using the fact that the suspension
of a PL-sphere is again a PL-sphere of one greater dimension).\hfill {\it q.e.d.}

\vskip .2in\leftskip .2in

\noindent {\bf Example 2}~If $g=0$, $r=s=1$, and $\delta _1=2$, then $F$ is a
once-punctured bigon.  Any arc in $F$ with distinct endpoints must be inessential, and any
arc whose endpoints coincide must be separating.  There are two
complementary components to $a$, one of which has one cusp in its
frontier, and the other of which has two cusps.  The component of the former type must
contain the puncture and of the latter type must be a triangle.  It follows that
$Arc(F)=Arc'(F)$ is in this case again a zero-dimensional sphere, and Lemma 13 applies
as before to prove sphericity whenever $g=0$ and $r=s=1$.

\vskip .2in\leftskip=0ex

\noindent 
As in this example, if $a$ is a separating arc whose endpoints coincide with the point $p$ in any surface,
then its two complementary components are distinguished by the fact that one component has exactly one 
cusp arising from $p$, and the other component has exactly two cusps arising from $p$, and we
shall refer to these respective complementary regions as the ``one-cusped'' and ``two-cusped'' components.   In
particular, any arc in a planar surface whose endpoints coincide is necessarily separating.

\vskip .2in\leftskip .2in

\noindent {\bf Example 3}~If $g=s=0$, $r=2$, and $\delta _1=\delta _2=1$, then $F$ is an
annulus with one distinguished point on each boundary component.  It is elementary that
any essential arc must have distinct endpoints and its isotopy class is classified by an integral ``twisting
number'', namely, the number of times it twists around the core of the annulus;
furthermore, two non-parallel essential arcs can be disjointly embedded in $F$ if and only
if their twisting numbers differ by one.  Thus, $Arc'(F)$ is in this case isomorphic to the
real line, the mapping class group $PMC(F)$ is infinite cyclic and generated by the Dehn twist along
the core of the annulus, which acts by translation by one on the real line, and the
quotient is a circle.  Again, Lemma~13 extends this to any surface with $g=s=0$ and $r=2$.

\vskip .2in

\noindent {\bf Example 4}~If $g=0$, $s=2$, and $r=1=\delta _1$, then $F$ is a
twice-punctured monogon.  If $a$ is an essential arc in $F$, then its one-cusped component
contains a single puncture, which may be either of the two punctures in $F$, and the
two-cusped component contains the other puncture of $F$.  Labeling the punctures of $F$
by 0 and 1, choose disjointly embedded arcs $a_i$ so that $a_i$ contains puncture $i$ in
its one-cusped complementary component, for $i=0,1$.  The Dehn twist $\tau$ along the
boundary generates $PMC(F)$, and the 0-skeleton of $Arc'(F)$ can
again be identified with the integers, where the vertex $\sigma (\{ \tau ^j(a_i)\} )$ is identified with
the integer $2j+i$~~for any integer $j$ and $i=0,1$.
 
\vskip .1in

\noindent Applying Example~2 to the
two-cusped component of $F_{\{a_i\}}$, for $i=0,1$, we find exactly two
one-simplices incident on
$\sigma (\{ a_i\} )$, namely, there are two ways to add an arc disjoint from $a_i$
corresponding to the two cusps in the two-cusped component, and these arcs are identified
with the integers $i\pm 1$.  Thus, $Arc'(F)$ is again identified with the real line,
$PMC(F)$ acts this time as translation by two, and $Arc(F)$ is again a circle.  Lemma~13
then proves sphericity of $Arc(F)$ for any planar surface with $r=1$ and $s=2$ 

\vskip .2in\leftskip=0ex

\noindent The typical aspect illustrated in this example is that the
arcs in an arc family come in a natural linear ordering which is invariant under the action
of the pure mapping class group.  Specifically, enumerate once and for all the
distinguished points $p_1,p_2,\dots ,p_\Delta$ in the boundary of $F$ in any manner. 
The boundary of a regular neighborhood of 
$p_i$ in
$F$, for $i=1,2,\dots , \Delta$ contains an arc $A_i$ in the interior of $F$, and $A_i$ comes equipped with a
natural orientation (lying in the boundary of the component
of
$F-A_i$ which contains $p_i$).
If
$\alpha $ is an arc family in $F$, then there is a first intersection of a
component of $\alpha $ with the $A_i$ of least index, and this is the first arc $a_1$ in the
linear ordering.  By induction, the $(j+1)^{\rm st}$ arc (if such there be) in the putative linear
ordering is the first intersection (if any) of a component of $\alpha -\{ 
a_1,\dots ,a_{j-1}\}$ with the $A_i$ of least index. 

\vskip .1in

\noindent Since we take the pure mapping class group which fixes each distinguished point
and preserves the orientation of $F$, this linear ordering descends to a well-defined
linear ordering on the $PMC(F)$-orbits in $Arc'(F)$.  It follows immediately that the
quotient map $Arc'(F)\to Arc(F)$ is injective on the interior of any simplex in $Arc'(F)$, so that
the simplicial complex $Arc'(F)$ descends to a CW decomposition on $Arc(F)$. 

\vskip .1in 

\noindent Furthermore, in the
action of
$PMC(F)$ on
$Arc'(F)$, there can be no finite isotropy, and the isotropy subgroup of a simplex $\sigma (\alpha )$
in
$PMC(F)$ is either trivial or infinite.
The former case occurs if and only if 
$\alpha $ {quasi fills} the surface $F$.

\vskip .1in

\noindent One fundamental inductive tool is:

\vskip .2in

\noindent {\bf Proposition 14} [P3;Lemma 3]\it Fix a bordered surface $F$, and suppose that $Arc(F_\alpha )$ is PL-homeomorphic to a
sphere of dimension
$N(F _\alpha)$ for each non-empty arc family
$\alpha $ in $F$.  Then $Arc(F)$ is an $N(F)$-dimensional PL-manifold.\rm

\vskip .2in

\noindent  \noindent{\bf Idea of Proof}~~ The idea is that the link of a simplex $\sigma [\alpha ]$ in $Arc(F)$ in the first
barycentric subdivision is canonically isomorphic to $Arc(F_\alpha )$; this is a general fact about geometric realizations of partially
ordered sets.  Thus, any point
$x\in Arc(F)$ admits a neighborhood which is PL-homeomorphic to an open ball of dimension $N=N(F)$.  Of course, $x$ lies in the interior
$Int~\sigma [\alpha ]$ of $\sigma [\alpha ]$, for some arc family $\alpha$ in $F$, and we suppose that $\alpha$ is comprised
of
$p+1$ component arcs. We have the identity $N+1= (N(F_\alpha )+1)+(p+1)$ since both sides of the equation give the
number of arcs in a quasi-triangulation of $F$.  Since $Arc(F_\alpha )$ is PL-homeomorphic to an $N(F_\alpha
)$-dimensional sphere by hypothesis, $Int~ \sigma [\alpha ]*Arc(F_\alpha )$ is PL-homeomorphic to an open ball of
dimension
$p+N(F_\alpha )+1 =N$.  This gives the required neighborhood of $x$ in $Arc(F)$ and proves that $Arc(F)$ is inductively a manifold.
~~~~~~~\hfill{\it q.e.d.}

\vskip .2in

\noindent  As was mentioned before, our basic approach to sphericity is inductive, and the proof depends upon the inductive
hypothesis that the sphericity conjecture holds for $F_\alpha$ for any non-trivial arc family $\alpha$ in $F$, i.e., by induction on
$N(F)$.  In order to avoid the Poincar\'e Conjecture in dimensions less than four, our basis step involves an explicit analysis of all
arc complexes of dimension at most four.  

\vskip .1in

\noindent In the special case $g=0=r-1$ of multiply punctured polygons, we have
already analyzed in the examples above (in combination with Lemma 13) all of the arc complexes of dimension at most two.  In this case, 
the only remaining complex in dimension three is for the surface $F_{0,1}^3$ and in dimension four is for $F_{0,(2)}^3$, which follows
from $F=F_{0,1}^3$ by Lemma 13.

\vskip .1in

\noindent To handle the arc complexes of $F=F_{0,1}^s$, for $s\geq 3$, label the various punctures of $F$ with distinct members of the set
$S=\{ 1,2,\dots ,s\}$.  If $\alpha $ is an arc family in $F$ and $a$ is a component arc of $\alpha $,
then $a$ has a corresponding one-cusped component containing some collection of punctures
labeled by a proper non-empty subset
$S(\alpha )\subseteq S$,
which we regard as the ``label'' of $a$ itself. 

\vskip .1in

\noindent More generally, define a ``tableaux'' $\tau$ labeled by $S$ to be a rooted tree embedded in
the plane where: the (not necessarily univalent) root of $\tau$ is an unlabeled vertex, and the other
vertices of $\tau$ are labeled by proper non-empty subsets of
$S$; for any $n\geq 1$, the vertices of $\tau$ at distance $n$ from the root are pairwise disjoint  subsets
of $S$; and, if a simple path in $\tau$ from the root passes consecutively through
the vertices labeled
$S_1$ and $S_2$, then $S_2$ is a proper subset of $S_1$.

\vskip .1in

\noindent Given an arc family $\alpha $ in $F$, inductively define the corresponding
tableaux
$\tau =\tau (\alpha )$ in the natural way:  For the basis step, choose as root some point in the component of
$F_{\alpha }$ which contains the boundary of $F$; for the inductive step, given a vertex of $\tau$
lying in a complementary region $R$ of $F_{\alpha }$, enumerate the component arcs
$a_0,a_1,\dots ,a_m$ of $\alpha $ in the frontier of $R$, where we assume that these
arcs occur in this order in the canonical linear ordering 
described in $\S$2 and $a_0$ separates $R$ from the root.  Each arc $a_i$ separates $R$ from another component
$R_i$ of
$F_{\alpha }$, and we adjoin to $\tau$ one vertex in each such component $R_i$ with the label $S(a_i)$ together
with a one-simplex connecting $R$ to $R_i$, for each
$i=1,2,\dots ,m$, where the one-simplices are disjointly embedded in $F$.  
\vskip .1in

\noindent It follows immediately from the topological classification of surfaces that
$PMC(F)$-orbits of cells in $Arc(F)$ are in one-to-one
correspondence with isomorphism class of tableaux labeled by $S$.
Furthermore, since the edges of $\tau (\alpha )$ are in one-to-one correspondence with the component
arcs of
$\alpha $, a point in $Int(\sigma (\alpha ))$ is uniquely
determined by a projective positive weight on the edges of $\tau(\alpha )$.  It follows that $Arc(F)$
itself is identified with the collection of all such projective weightings on all isomorphism classes of
tableaux labeled by $S$. 

~~\vskip 3.8in

\hskip .3in {{{\epsffile{f0131.epsf}}}}

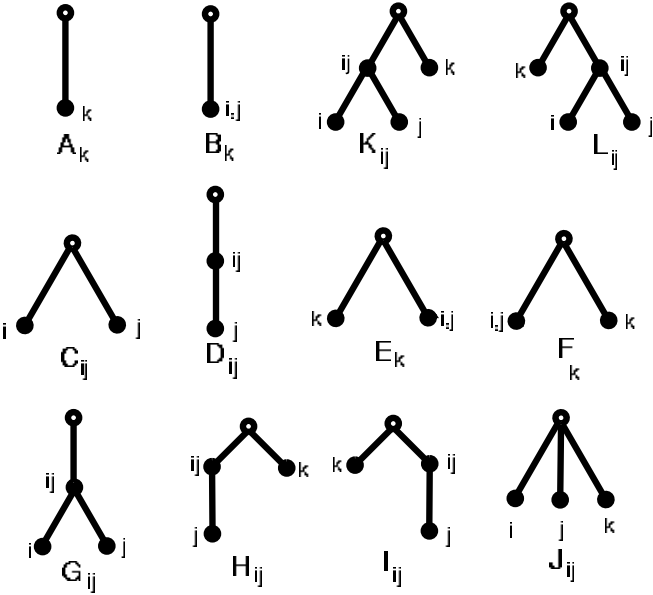
\centerline{{\bf Figure 7}~~{The tableaux for $F_{0,1}^3$.}}

\vskip .2in

\leftskip .2in

\noindent {\bf Example 5}~~The surface $F_{0,1}^3$.  Let us adopt the convention that given an ordered pair $ij$, where $i,j\in\{
1,2,3\}$, we shall let $k=k(i,j)=\{ 1,2,3\} -\{ i,j\}$, so $k$ actually depends only upon the unordered pair
$i,j$.  The various tableaux for $F_{0,1}^3$ are enumerated, labeled, and indexed in Figure~7, where in each case, $ij$ varies
over all ordered pairs of distinct members of $\{ 1,2,3\}$, $k=1,2,3$, and the bullet represents the root.  
In this notation and
letting
$\partial$ denote the boundary mapping in $Arc=Arc(F_{0,1}^3)$, one may directly compute incidences of cells summarized as follows.

\settabs 5\columns

\vskip .2in

\+&$\partial C_{ij}=A_i-A_j$,&&$\partial H_{ij}=D_{ij}-C_{jk}+F_k$,\cr
\+&$\partial D_{ij}=A_j-B_k$,&&$\partial I_{ij}=D_{ij}-E_k+C_{kj}$,\cr
\+&$\partial E_{k}=A_k-B_k$,&&$\partial J_{ij}=C_{ij}-C_{ik}+C_{jk}$,\cr
\+&$\partial F_{k}=B_k-A_k$,&&$\partial K_{ij}=G_{ij}-H_{ij}+H_{ji}-J_{ij}$,\cr
\+&$\partial G_{ij}=C_{ij}-D_{ji}+D_{ij}$,&&$\partial
L_{ij}=I_{ij}-G_{ij}+J_{ki}-I_{ji}$.\cr

\vskip .2in

\noindent We may symmetrize and define sub complexes ${X}_{ k}=X_{ij}\cup X_{ji}$, for $X=K,L$, and furthermore set
${M}_{ k}={ K}_{k }\cup {L}_{k}$, f
or $k=1,2,3$.  Inspection of the incidences of cells shows that each of ${K}_k$ and
${L}_k$ is a 3-dimensional ball embedded in $Arc$, as illustrated in Figures~8a and 8b, respectively,
with $K_{ij}$ on the top in part a) and $L_{ji}$ on the top in part b).  Gluing together $L_k$ and $K_k$ along their common
faces $G_{ij}, G_{ji}$, we discover that
${M}_{k}$ is almost a 3-dimensional ball embedded in
$Arc$ except that two points in its boundary are identified to the single point $A_k$ in $Arc$, as illustrated in Figure~8c.

~~\vskip 1.8in

~{{{\epsffile{f0132.epsf}}}}

\vskip .1in

\noindent \hskip .5in Figure~8a ${K}_{k}$.\hskip .6in Figure 8b ${ L}_{k}$.\hskip .5in Figure 8c ${ M}_{k}$.

\vskip .1in

\centerline{{\bf Figure 8}~~{${ M}_{k}$ and the balls ${K}_{k},{ L}_{k}$.}}

\vskip .2in

\noindent Each ${ M}_k$, for $k=1,2,3$, has its boundary entirely contained in the sub complex
${J}$ of $Arc$ spanned by $$\{ A_k: k=1,2,3\}\cup\{ C_{ij},J_{ij}: i,j\in\{ 1,2,3\}~{\rm are~distinct}\}.$$  In order to understand
${J}$, we again symmetrize and define $J_{k}=J_{ij}\cup J_{ji}$, so each $J_{k}$ is isomorphic to a cone, as
illustrated in Figure 9; 
we shall refer to the 1-dimensional simplices $C_{ik},C_{jk}$ as the ``generators'' and to
$C_{ij},C_{ji}$ as the ``lips'' of
$J_{k}$.  The one-skeleton of ${ J}$ plus the cone $J_{k}$ is illustrated in Figure 9.
Imagine taking $k=3$ in Figure 9 and adjoining the cone $J_2$ so that  
the generator $C_{12}$ of $J_{2}$ is attached to the lip $C_{12}$ of $J_{3}$, and the generator $C_{13}$ of
$J_{3}$ is attached to the lip $C_{13}$ of $J_{2}$.  Finally, ${ J}$ itself is produced by symmetrically attaching
$J_{1}$ to $J_{2}\cup J_{3}$ in this lip-to-generator fashion.

~~\vskip 1.5in

\hskip .8in{{{\epsffile{f0133.epsf}}}}

\vskip .1in

\centerline{{\bf Figure 9}~~{The cone ${ J}_{k}$ and the one-skeleton of ${ J}$.}}

\vskip .2in

\noindent In order to finally recognize the 3-dimensional sphere, it is best to take a regular neighborhood of ${ J}$ in
$Arc$, whose complement is a disjoint union of three 3-dimensional balls.  Each 3-dimensional ball is naturally identified with
the standard ``truncated'' 3-simplex, where a polyhedral neighborhood of the 1-skeleton of the standard 3-simplex has been excised. 
These truncated simplices are identified pairwise along pairs of hexagonal faces to produce the 3-dimensional sphere in the natural way.

\leftskip=0ex

\vskip .2in

\noindent We next prove that the arc complex $Arc(F)$ minus a point is contractible, for any surface
$F=F_{0,(1)}^s$.  Let $a$ denote an arc labeled by $\{
s\}$, and let
$b$ denote an arc labeled by $S-\{ s\}$.  We shall show that ${\cal A}=Arc(F)- \{ [a]\}$ strong
deformation retracts onto $ \{ [b]\}$.  The retraction proceeds in four steps, as follows: ${\cal A}$
retracts to
$${\cal B}=\bigl \{ [\alpha ]\in{\cal A}: ~{\rm no~vertex~of}~\tau (\alpha )~{\rm
is~labeled~by}~\{ s\}\bigr\} ,$$
${\cal B}$ retracts to
$${\cal C}=\bigl \{ [\alpha ]\in{\cal B}: s~{\rm
is~not~a~member~of~the~label~of~any~vertex~of}~\tau (\alpha )\bigr\} ,$$
${\cal C}$ retracts to
$${\cal D}=\{ [\alpha ]\in{\cal C}: \exists !~{\rm
vertex~of}~\tau (\alpha )~{\rm adjacent~to~the~root~labeled~by}~S-\{ s\}\bigr\},$$
and finally ${\cal D}$ retracts to the point $[b]$.

\vskip .1in

\noindent In order to coherently describe these retractions, let us introduce a natural linear
ordering, the so-called ``pre-ordering'', on the vertices of any tableaux $\tau$.  Begin
with the root of $\tau$ as the first element in the putative pre-ordering.  Given a vertex $v$ at distance
$n\geq 0$ from the root, there is a (possibly empty) family
$v_1,v_2,\dots ,v_m$ of vertices of $\tau$ at distance $n+1$ from the root which are adjacent to $v$, and this
family of vertices comes in the right-to-left linear ordering from the definition of a tableaux. 
Furthermore, each $v_j$ is the root of a (possibly empty) sub tableaux $\tau _j$ of
$\tau$, for $j=1,2,\dots ,m$, whose other vertices (if any) are at distance at least $n+2$ from the root in
$\tau$.  Take the pre-ordering on $\tau _1$, then on $\tau _2$, and so on up to $\tau _m$, to
complete the recursive definition of the pre-ordering on the vertices of $\tau$ itself.

\vskip .1in

\noindent For the first retraction from ${\cal A}$ to ${\cal B}$, 
simply decrease the projective weights, one at a time, on the one-simplices of $\tau (\alpha )$, for $[\alpha
]\in {\cal A}-{\cal B}$, which terminate at vertices labeled by $\{ s\}$ beginning with the greatest such vertex
in the pre-ordering.

\vskip .1in

\noindent For the second retraction from
${\cal B}$ to ${\cal C}$, consider the greatest vertex $v$ of $\tau (\alpha )$, for $[\alpha ]\in{\cal B}-{\cal
C}$, in the pre-ordering whose label
$S$ contains $s$, and let $e$ be the edge of $\tau (\alpha )$ terminating at $v$, where $e$ has a
projective weight $w$.  Insert a new bi-valent vertex $u$ into $e$ labeled by $S$ and re-label $v$ by
$S-\{ s\}$.  The edge $e$ is thus decomposed into two edges, and the one incident on $v$ is given
projective weight $tw$, while the other is given projective weight $(1-t)w$, for $0\leq t\leq 1$. 
The corresponding deformation as $t$ goes from 0 to 1 effectively replaces the label $S$ on $v$ with
the label $S-\{ s\}$ and thus reduces the number of vertices whose label contains $s$.  Iteratively
removing $s$ from the label of the greatest vertex in this manner defines the retraction of ${\cal B}$
to ${\cal C}$.

\vskip .1in

\noindent For the third retraction, suppose that
$[\alpha ]$ lies in ${\cal C}-{\cal D}$, and let $v$ denote the root of $\tau (\alpha )$.  Adjoin to
$\tau (\alpha )$ a new edge
$e$ incident on $v$, and let the other endpoint of $e$ be the new root of the resulting tableaux,
where 
$v$ is labeled by $S-\{ s\}$.  Let $w$ denote the projective class of the sum of all of the weights
on $\tau (\alpha )$, and define a projective weight $tw$ on $e$, for $0\leq t\leq 1$.
The corresponding deformation as $t$ goes from 0 to 1 effectively adds a new vertex with label
$S-\{ s\}$.

\vskip .1in

\noindent Finally, simply decrease the projective weights on all of the edges of $\tau (\alpha )$
except the edge incident on the root in order to contract ${\cal D}$ to $[b]$, completing the proof that the complement of a point in
the arc complex $Arc(F_{0,(1)}^s)$is contractible.

\vskip .1in

\noindent

\noindent We shall say that a surface $F$ ``satisfies the inductive hypothesis'' provided
$Arc(F_\alpha )$ is PL-homeomorphic to a
sphere of dimension
$N(F _\alpha)$ for each non-empty arc family
$\alpha $ in $F$.  (For the perspicacious reader, notice that we tacitly assumed this inductive hypothesis in the proof of Lemma 13.) 
The inductive hypothesis thus implies that
$Arc(F)$ is a manifold by Proposition~14, and this manifold is simply connected and a homology sphere since the complement of a point
is contractible by the previous discussion.  By the Poincar\'e Conjecture in high dimensions, it follows inductively that
$Arc(F)$ is a sphere for $F=F_{0,(1)}^s$, and hence for $F=F_{0,(n)}^s$ for any $n\geq 1$ by Lemma 13.  

\vskip .1in
 \noindent This completes the promised
proof of the sphericity theorem for punctured polygons.

\vskip .3in

\noindent {\bf 6. Punctured surfaces and fatgraphs}\vskip .2in

\noindent We take this opportunity to point out that in [P3] we claimed that our
compactification maps continuously to the Deligne-Mumford compactification, and this may not be true; the precise connection between
the two compactifications will be taken up elsewhere.  Secondly, Theorem 5 of [P3] is slightly corrected by Theorem 17 below.

\vskip .1in

\noindent We return to the setting $F=F_{g,0}^s$ of punctured surfaces with $s\geq 1$ and $r=0$.  Let us choose from among the punctures
of
$F$ a distinguished one, to be denoted $\pi$.  Build as before the geometric realization ${\cal A}_\pi '(F)$ of the partially ordered
set of all isotopy classes of arc families of simple essential arcs which are asymptotic to $\pi$ in both directions (referred to
subsequently as arcs ``based'' at $\pi$).  Consider the punctured {\it arc complex of $F=F_g^s$} given by ${\cal A}_\pi (F)={\cal A}_\pi
'(F)/PMC(F)$.  Notice that given an arc family $\alpha$ in $F$, each component of $F_\alpha$ inherits a labeling as a bordered surface,
where the ends of $F_\alpha$ are taken to be the distinguished points in the boundary; this is the basic connection between bordered and
punctured arc complexes.

\vskip .1in

\noindent Several of the main results of [P3] are summarized in the following theorem.

\vskip .2in

\noindent{\bf Theorem 15}~~\it For any choice of puncture $\pi$, the Teichm\"uller space $T_g^s$ of $F=F_g^s$ is naturally isomorphic
to the subspace of
${\cal A}'_\pi={\cal A}'_\pi(F)$ corresponding to arc families based at $\pi$ which quasi-fill $F$, so ${\cal A}_\pi$ is a
compactification of $M=M_g^s$.  The third barycentric subdivision of ${\cal A}'_\pi$ descends to an honest simplicial complex on the
compactification ${\cal A}_\pi$ of $M$. \rm

\vskip .2in

\noindent Degenerate structures (i.e., ``ideal'' points of ${\cal A}_\pi-M$) will be fully described presently as will a graphical
description of cells together with the corresponding
matrix-model.  For the moment, though, we turn to sphericity.  In fact, orbifoldicity of ${\cal A}_\pi$ follows from the
sphericity conjecture.

\vskip .1in

\noindent Indeed, as in the proof of Proposition~14, where the link of a simplex $\sigma [\alpha ]$ in the second barycentric
subdivision of
$Arc(F)$ is isomorphic to $Arc(F_\alpha )$, this time there is a finite group ${\cal G}_\alpha$ acting on $\sigma [\alpha ]*Arc
'(F_\alpha )$, and the natural mapping $Int\bigl (\sigma [\alpha ]*Arc'(F\alpha )\bigr )\to {\cal A}'(F)$ induces a homeomorphism
between its image and the quotient $Int\bigl (\sigma [\alpha ]*Arc'(F\alpha )\bigr )/{\cal G}_\alpha$.  Thus, the links of simplices
are virtually arc complexes of bordered surfaces (i.e., spheres), so ${\cal A}_\pi$ is indeed an
orbifold.  The local groups
${\cal G}_\alpha$ in the orbifold structure of ${\cal A}_\pi$ are easily described and recognized as follows.  Consider the stabilizer
$S_\alpha <PMC(F)$ of
$F_\alpha$ in
$F$ (which may permute the labeling of $F_\alpha$ as a bordered surface), and notice that the quotient
$${\cal G}_\alpha =S_\alpha /PMC(F)$$
is a finite group.  For instance, a component of $F_\alpha$ might be an annulus of type $F_{0,(1,n)}^0$, and there may be some $\phi\in
PMC(F)$ which cyclically permutes the $n$ distinguished points on a common boundary component;
the resulting coset corresponds to a ``fractional Dehn twist''. 

\vskip .1in

\noindent The starting point for the sphericity conjecture was the request/challenge from Ed Witten in 1989 to find an
orbifold compactification of moduli space together with a cell decomposition which restricts to the usual one on moduli space itself.
The sphericity theorem for multiply punctured disks (proved in $\S$5) already proves orbifoldicity for punctured spheres as
noted in [P3].  The sphericity conjecture then gives orbifoldicity of our new compactification for any punctured surface.
\vskip .1in

\noindent Given a quasi triangulation $\tau$ of $F=F_g^s$ based at $\pi$, we may formally add in a choice-free manner an edge to
$\tau$ in each punctured component of $F_\tau$ to produce an ideal triangulation $\tau ^+$ of $F$.  
There is a fatgraph $G$ in $F$ which is Poincar\'e dual to $\tau ^+$ in $F$.  On this level of dual fatgraphs, there is one loop of $G$
for each puncture of $F$ other than $\pi$; these loops are taken formally here in the sense that we associate a lambda length to each
such loop, where these lambda lengths are taken to be large relative to the non-loop lambda lengths and are otherwise unrestricted.

\vskip .1in

\noindent The promised description of degenerate structures is given by the next theorem.

\vskip .2in

\noindent {\bf Theorem 16} [P3; Theorem 7]~~\it Fix a quasi triangulation $\tau$ of $F_g^s$ with associated ideal triangulation $\tau
^+$.  Let
$E_j(t)$ be a degeneration of simplicial coordinates on $\tau ^+$, and let $\beta$ denote the collection of arcs in $\tau ^+$ so that
the corresponding lambda lengths $e_i$ do not diverge.  Fixing a horocycle about $\pi$, each component of $F_\beta$ comes equipped with
a canonical decoration at the cusps on the boundary.  The polygons and once-punctured polygons of $F_\beta$ combine to give a decorated
hyperbolic structure on a bordered subsurface $GF$, and the complement $TF=F-GF$ is regarded as a topological bordered surface
(with no hyperbolic structure) together with a decoration coming from the specified horocycle.\rm

\vskip .2in

\noindent Thus, a point in ${\cal A}_\pi$ gives rise to a geometric part $GF\subseteq F$ together with a topological part $TF\subseteq F$,
and the topological part is non-empty if and only if the point is ideal.  One is thus permitted to ``forget'' the hyperbolic structure on
a subsurface $TF$.  In particular, the Deligne-Mumford compactification corresponds to ``forgetting'' the structure only on disjoint
collections of annuli $TF\subseteq F$.

\vskip .1in

\noindent An elementary observation is that any graph $G$ admits a canonical decomposition into its recurrent part $RG$ and its
non-recurrent part $NG$, as follows.  We shall say that a closed edge-path $P=e_1,e_2,\ldots ,e_n$ on $G$, thought of as a
concatenation of edges, is ``edge-simple'' provided:

\vskip .1in

\leftskip .3in

\noindent ~i) $e_i\neq e_{i+1}$ for any $i$ (where the index $i$ is taken modulo $n+1$);

\vskip .1in

\noindent ii) $P$ never twice traverses the same {\sl oriented} edge of $G$.

\leftskip=0ex

\vskip .1in

\noindent The decomposition is given by:
$$\eqalign{
RG&=\{e\in G:\exists~{\rm edge-simple~curve~through}~e\},\cr
NG&=G-RG.\cr
}$$

\vskip .1in

\noindent The next result improves and slightly corrects an earlier result.

\vskip .2in

\noindent {\bf Theorem 17} [P3; Theorem 5]~~\it Consider a degeneration of hyperbolic structure, where simplicial coordinates $E_i(t)\to
0$ for the indices
$i\in I$ and lambda lengths $e_j(t)\to\infty$ for indices $j\in J$.  Then $J\subseteq I$.  Furthermore, letting $G_I=\{ e\in G:e\in
I\}$ and likewise $G_J$, we have $R(G_I)=G_J$.\rm

\vskip .2in

\noindent Indeed, that $J\subseteq I$ follows directly from the second part of Lemma~5, and the first part of Lemma~5 is applied to
prove the remaining assertions about recurrent subgraphs.

\vskip .2in

\noindent  The combinatorics of our compactification is clear on the level of arc families based at $\pi$ in a surface $F=F_g^s$.  Cells
simply correspond to arc families which may or may not quasi fill, and the topological part of the decomposition of $F$, if any,
arises from the complementary components which are different from (perhaps once-punctured) polygons.  The graphical formulation is
somewhat complicated but interesting nonetheless, and we finally describe it.  Given the quasi triangulation $\tau$ of $F$, complete
it to an ideal triangulation $\tau ^+$ as before with corresponding dual marked fatgraph $G'=G'_\tau$.  In contrast to the previous
discussion, we here simply remove all the loops to produce a marked fatgraph $G=G_\tau$, where $G$ has $s-1$ external nodes, and we
label these external nodes with the indices $2,\ldots ,s$ of the corresponding punctures.

\vskip .1in

\noindent One considers Whitehead moves (i.e., contraction/expansion of edges of the fatgraphs) as usual, and now also another
elementary move supported on a surface of type $F_{0,(2)}^2$.  This new move interchanges the two possible quasi triangulations of this
surface, so on the level of dual fatgraphs, given a puncture $i$, the new move just flips the edge containing the terminal node labeled
$i$ by altering the fattening at its other vertex.  In fact, this new elementary move is just the composition of two Whitehead moves on
ideal triangulations of $F_{0,(2)}^2$.

\vskip .1in

\noindent We now two-color the edges of $G$ by specifying ``regular'' and ``ghost'' edges, where we thing of ghost edges as ``missing'',
and there must be at least one regular edge.

\vskip .1in

\noindent Take equivalence classes of these two-colored fatgraphs with labeled terminal nodes under Whitehead moves and the new
elementary move, where we are only allowed to perform moves along the ghost edges.  Furthermore, once an edge is a ghost edge, it remains
a ghost (``once you're dead, you're dead''), and otherwise the color of an edge is unaffected.  For instance, in order to
collapse a single edge with distinct vertices in a cubic fatgraph to a quartic vertex, instead, put a ghost icon on the edge to be
collapsed; collapse and expand it to produce another cubic fatgraph with a ghost icon on the new edge, and identify these two ghostly
fatgraphs.

\vskip .2in

\noindent {\bf Theorem 18} [P3; Theorem 9]~~\it The dual of the simplicial complex ${\cal A}'_\pi$ is naturally isomorphic to the
geometric realization of these marked two-colored labeled fatgraph equivalence classes.\rm

\vskip .2in

\noindent An obvious project is to compute the virtual Euler characteristics
using
this two-colored matrix model.  It is also worth mentioning that the graphical formulation of $Arc(F)$ for $F$ a bordered surface can be
described in analogy to this discussion, where each distinguished point in the boundary gives rise to a ``hair'' in the natural way.  The
dual of $Arc'(F)$ is naturally isomorphic to the set of two-colored fatgraphs with $s$ external nodes and $\Delta$ hairs with external
nodes, both labeled, modulo Whitehead moves on the ghost edges, where each boundary
component has at least one hair and hairs are never allowed to be ghosts (which has a macabre corporeal interpretation).

\vskip .3in

\noindent {\bf 7. Operads}\vskip .2in

\noindent Several topological and homological operads
based on families of projectively weighted arcs in bordered surfaces $F_{g,r}^s$ are
introduced and studied in [KLP].  This work as well as further material [K1] due to Kaufmann
are surveyed in his paper [K2] in this volume, so our treatment here will be brief and contextual.  The spaces
underlying the basic operad are identified with open subsets of $Arc(F)$ which contain $Arc_\# (F)$ for suitable
families of bordered  surfaces $F$.  Specifically, say that an arc family $\alpha$ in $F$ is {\it exhaustive} if  
for every boundary component of $F$, there is some arc in $\alpha$ which is incident on this boundary component, and
define
$$
Arc_g^s(n)=\{ [\alpha ]\in
Arc(F_{g,n+1}^s): \alpha ~{\rm is~exhaustive}\} .
$$
We shall see that ``gluing in the spirit of train tracks'' imbues $Arc=\{Arc_g^s(n)\}_{n\geq 1}$ with a cyclic
operadic structure which restricts to a sub-operadic structure on $Arc_{cp}=\{Arc_0^0(n)\}_{n\geq 1}$; there are also
natural twistings $\widetilde{Arc}_g^s(n)=Arc_g^s(n)\times (S^1)^{n+1}$ with related operadic structures for which we
refer the reader to [KLP]; this trivial Cartesian product
is interpreted as a family of coordinatizations of
the boundary of $F$ which are twisted by $Arc_g^s(n)$.

\vskip .1in

\noindent Already algebras
over $Arc_{cp} < Arc$ are shown in [KLP] to be Batalin-Vilkovisky algebras where the
entire BV structure is realized simplicially.
An example of this simplicial realization is given by the symmetric diagram which describes
the BV equation itself, as described in [K2].  Furthermore, this basic operad contains the
Voronov cacti operad [Vo] up to homotopy.

\vskip .1in

\noindent Our operad composition on $Arc$ depends upon
an explicit method of
combining families of {\sl projectively weighted} arcs in surfaces, that is, each component arc of the family is
assigned the projective class of a positive real number; we
next briefly describe this composition.  Suppose that
$F^1,F^2$ are surfaces with distinguished boundary components $\partial ^1\subseteq F^1, \partial
^2\subseteq F^2$.  Suppose further that each surface comes equipped with a properly embedded family of
arcs, and let
$a_1^i,\ldots ,a_{p^i}^i$ denote the arcs in
$F^i$ which are incident on $\partial ^i$, for $i=1,2$, as illustrated in part I of Figure~10.
Identify $\partial ^1$ with $\partial ^2$ to produce a surface $F$.  We wish to furthermore combine
the arc families in $F^1,F^2$ to
produce a corresponding arc family in $F$, and there is evidently no well-defined way to achieve this
without making further choices or imposing further conditions on the arc families (such as $p^1=p^2$).

\vskip .1in

\noindent The additional data required for gluing is given by an assignment of the projective class of one real
number, a weight, to each arc in each of the arc families.  
The weight $w_j^i$ on $a_j^i$ is interpreted geometrically as the height of a rectangular band
$R_j^i=[0,1]\times [-{{w_j^i}/ 2}, {{w_j^i}/ 2}]$ whose core $[0,1]\times\{ 0\}$ is identified with $a_j^i$, for $i=1,2$ and
$j=1,\ldots ,p^i$. If we assume that
$\sum _{j=1}^{p^1} w^1_j=\sum _{j=1}^{p^2} w^2_j$ for simplicity, so that the total height of all the bands
incident on $\partial ^1$ agrees with that of $\partial ^2$; in light of this assumption, the bands in
$F^1$ can be sensibly attached along $\partial ^1$ to the bands in $F^2$ along $\partial ^2$ to
produce a collection of bands in the surface $F$ as illustrated in part II of Figure~10;
notice that the horizontal edges of the rectangles $\{ R_j^1\} _1^{p^1}$ decompose the rectangles 
$\{ R_j^2\} _1^{p^2}$ into sub-rectangles and conversely.  The resulting family of sub-bands, in turn,
determines a weighted family of arcs in $F$, one arc for each sub-band with a weight given by the width
of the sub-band; thus, the weighted arc family in $F$ so produced depends upon the weights in a
non-trivial but combinatorially explicit way.  This describes the basis of our gluing operation on families of weighted arcs, which is
derived from the theory of train tracks and partial measured foliations [PH].  In fact, the simplifying
assumption that the total heights agree is obviated by considering not weighted families of arcs, but
rather projectively weighted families of arcs since we may
de-projectivize in order to arrange that the simplifying assumption is in force (since $p^1p^2\neq
0$ by exhaustiveness), perform the construction just described, and finally re-projectivize.
The gluing operation just described may give rise to bands which are embedded annuli foliated by simple curves,
and these are simply discarded to finally produce a resulting projectively weighted arc family which is again
exhaustive.

~\vskip .2in

~~\vskip 2in

~{{{\epsffile{introfig.epsf}}}}

\vskip .1in

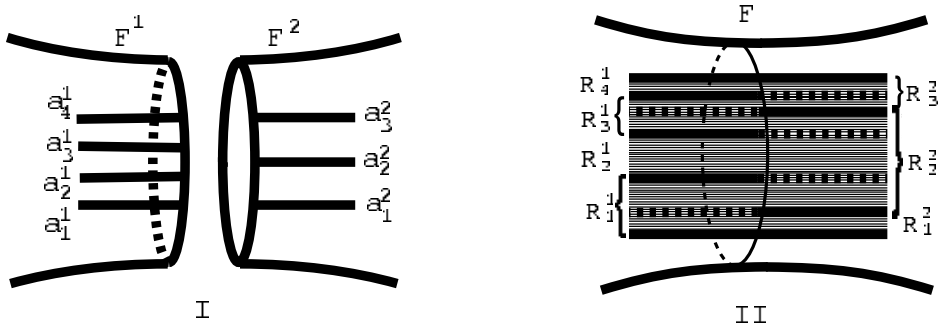
\centerline{{\bf Figure 10} Gluing weighted arc families as train tracks}

\vskip .2in
 
\noindent Enumerate the boundary components of $F^i$ as $\partial _0^i,
\partial _1^i,\ldots ,\partial _{n_i}^i$ once and for all for $i=1,2$; the gluing operation applied to boundary
components $\partial _i^1$ and $\partial _0^2$  induces the usual $\circ _i$ operation for $i=1,\ldots ,n_1$.
Letting
${\bf S}_p$
denote the $p^{\rm th}$ symmetric group, there is a natural
${\bf S}_{n+1}$-action on the labeling of boundary components which
restricts
to a natural
${\bf S}_n$-action on the boundary components labeled $\{ 1,2,\ldots ,
n\}$.
Thus, ${\bf S}_n$ and ${\bf S}_{n+1}$ act on $Arc_g^s(n)$, and extending by the
diagonal action of ${\bf S}_{n+1}$ on $(S^1)^{n+1}$, the symmetric groups
${\bf S}_n$
and ${\bf S}_{n+1}$ likewise act on $\widetilde{Arc}_g^s(n)$, where ${\bf S}_n$ by
definition acts trivially on the first coordinate in $(S^1)^{n+1}$.

\vskip .1in

\noindent {\bf Theorem 19} [KLP; Theorem 1.5.2-3]~~\it  The compositions $\circ _i$ are well-defined and imbue $Arc_{cp}$ with the
structure of a topological operad under the natural ${\bf S}_n$-action on labels.  The operad has a unit and is cyclic for the
natural ${\bf S}_{n+1}$-action.
$Arc$ likewise inherits the structure of a cyclic operad with unit.\rm

\vskip .2in

\noindent Just to give the flavor, we illustrate as one-parameter families of weighted arc families in Figure~11 the various operations
(the identity, BV operator, dot product, and star product), where the parameter $s$ satisfies $0\leq s\leq 1$ and the outermost
boundary component is labeled zero by convention.  In this same spirit, the BV equation itself is given by an explicit and symmetric
chain homotopy of weighted arc families as explained in [KLP] or [K2].

~~\vskip 1.1in

{{{\epsffile{ops.epsf}}}}

\vskip .1in

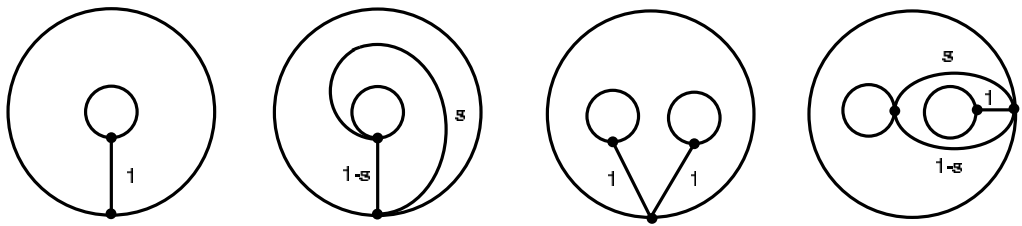
\centerline{{\bf Figure 11}~The identity, BV operator, dot product, and star operator}

\vskip .1in

\noindent Just as for cacti, $Arc_{cp}$ acts on the loop space of any manifold; does this action
extend to an action of $Arc$ or $\widetilde{Arc}$ itself?  Does our formalism allow for a coherent description and calculation of the
full string prop or its compactification?  Does the BV structure of $Arc$ or $Arc_{cp}$ survive injectively on the level of
homology?

\vskip .1in

\noindent As to the calculation of the homology of $Arc_{cp}$, the planar sphericity conjecture is of utility
since the complement of $Arc_{cp}(n)$ in $Arc(F_{0,n+1}^0)$ is a finite
union of spheres
$Arc(F_{0,r}^s)$, where $r+s=n+1$.  This arrangement of spheres should be analyzed using a Mayer-Vietoris
argument, and the homology of $Arc_{cp}$ itself could then be computed using Alexander duality in the sphere
$Arc(F_{0,n+1}^0)$; the twistings $\widetilde{Arc}_{cp}$ could also then be computed using the K\"unneth formula. 
Of course, the full sphericity theorem is likewise useful in
calculating the homology operad of
$Arc$ itself, of twistings such as
$\widetilde{Arc}$, and presumably also of the string prop.

\vskip .5in

\noindent {\bf Appendix. Biological Applications}\vskip .2in

\noindent The prediction of ``macromolecular folding'' is a central problem in contemporary biology.  The
``primary structure'' (i.e.,  the sequence of nucleotides $C,G,A,U$ for RNA, or of 20 amino acids for a general protein) can be
determined experimentally, and one wishes to predict from it the full three-dimensional structure or ``folding'' of the macromolecule. 
We have glossed over interesting and important mathematics here, the problem of ``sequence alignment'', one aspect of which is to
re-construct the entire primary structure of the macromolecule from empirically known snippets.

\vskip .1in

\noindent The basic model of the ``denatured
macromolecule'' is
$I=I_m=[0,m]\subseteq {\bf Z}\subseteq {\bf R}$ for some
$m$. 
A general ``bonding'' on $I_m$ is simply a collection ${\cal B}$ of unordered pairs $\{ i,j\}$, where $i,j\in{\bf Z}$, $|i-j|\geq 1$,
and $0\leq i,j\leq m$.  Various pairs of ``sites'' $i,j\in I$ are ``bonded'' if there is a corresponding element $\{ i,j\}\in{\cal
B}$, and we may consider the semi-circle $C_{ij}\subseteq {\bf R}\times {\bf R}_{>0}\subseteq {\bf R}^2$ with endpoints $(i,0),(j,0)$. 
One thinks of
${\bf R}$ as the ``backbone'' of the macromolecule, a tensile rod which bends to allow the bonding, where each $C_{ij}$ is collapsed to
a point; in particular, the backbone is sufficiently rigid to prohibit any bond $\{ i,j\}$ with $|i-j|=1$.  We say that a
bonding
${\cal B}$ is a ``secondary structure'' if
$C_{ij}\cap C_{k\ell}\subseteq {\bf R}$ for all
$\{ i,j\} ,\{ k,\ell \} \in{\cal B}$, whence $\{ C_{ij}:\{ i,j\} \in{\cal B}\}$ is an arc family in
${\bf R}\times{\bf R}_{\geq 0}\supseteq {\bf Z}\times {\bf R}_{\geq 0}\supseteq I_m\times {\bf R}_{\geq 0}$.  The secondary structure is
said to be ``binary'' provided $C_{ij}\cap C_{k\ell}\neq\emptyset$ implies $\{ i,j\}=\{ k,\ell \}$.
A ``helix'' of length $\ell$ is a collection $\{i,j\} ,\{ i+1,j-1\} ,\ldots \{ i+\ell ,i-\ell \}\in{\cal
B}$, whence it is a sequence of consecutive bonds.
There is a natural binary secondary structure associated with an arbitrary one as illustrated in Figure~12.

~~\vskip 1.1in

\leftskip=0ex

\hskip .4in{{{\epsffile{bnsffig2.epsf}}}}

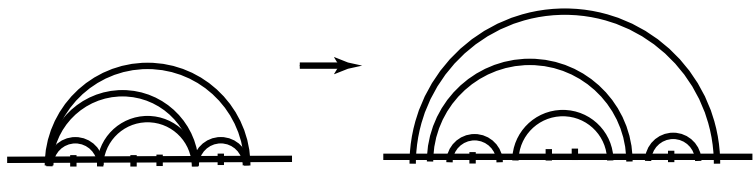
\centerline{{\bf Figure 12}~~{Binary from arbitrary foldings}}

\vskip .1in

\vskip .1in

\noindent An ``RNA secondary structure'' is a binary secondary structure on $I_m$ for some $m$.  This evidently depends
upon a plane of projection, but in practice, there seems to be a more or less well-defined plane of projection for the arc family.
(This is the plane the biologists draw for RNA which contains most of the bonding.)
In any case, the definition of an RNA secondary structure is mathematically complete.  In practice in ${\bf R}^3$, there is
further bonding between sites which does not lie in the plane of projection, and this is the
``tertiary structure'' or full three-dimensional folding, of which we shall give mathematical formulations later.  The prediction of
tertiary from primary structure for RNA is the full folding problem.  The prediction of the biologically correct secondary
structure from the primary structure for RNA is another example of the macromolecular folding problem, which we shall discuss
further below.

\vskip .1in

\noindent For a general protein, the ``secondary structure'' is a decomposition of it into two standard motifs: the $\alpha$-helix
and the $\beta$-sheet.  The $\alpha$-helix is an (essentially always right-handed) helix, and the $\beta$-sheet is a motif where the
macromolecule repeatedly folds back upon itself in a plane.  (There are other basic motifs for general proteins which can be added to
the definition of its secondary structure, but we shall not take this up here.)  Again in any case, the full three-dimensional
folding of a general protein is its ``tertiary stucture'', and the folding problem is to predict secondary or tertiary from
primary structure.

\vskip .1in

\noindent In the context of general proteins,
the most effective methods of folding prediction are based upon motif searches and matches with a database of primary structures whose
foldings have been experimentally determined by X-ray crystallography or other methods (and these experiments are delicate and
relatively few).  This is obviously disappointing to a theorist seeking {\sl ab initio} methods based on physical or mathematical
principles.  Any serious such model should evidently incorporate what seems to be the main driving force of
``hydophobicity/hydrophillicity'' (the tendency of individual nucleotides, amino acids or sequences of them to shun/seek water in the
biological environment), but experts seem to think that an effective thermodynamic model based on hydrophobicity is just
pie-in-the-sky.  Furthermore, structural considerations determined by the stoichiometry of the primary structure also seem to play a
role.  Effective prediction for biologically interesting primary sequences of length $10^2-10^9$ may be too complex
for any {\sl ab initio} methods.

\vskip .1in

\noindent For RNA secondary structures, however, there are rigorous and reasonably reliable
thermodynamic models, the ``standard'' one being Mfold [ZMT,ZS] which is defined as follows.  The input is the primary structure
(that is, a long sequence of letters $C,G,A,U$); the output is not only the optimal (i.e., global energy minimizing) folding or
foldings but also the ``sub-optimal'' band of the local minima whose values are within, say, five percent of the global minimum. 
Experience seems to indicate that the correct foldings at least lie in such a sub-optimal band.  The thermodynamics is governed by
``nearest neighbor bonds'', that is, Boltzmann weights are assigned to bonds and to pairs and triples of bonds in a helix at adjacent
sites; the energies of these configurations at fixed temperature have been determined experimentally, and the thermodynamic model can
be solved (for instance stochastically) on the computer.
This model apparently very accurately predicts melting points of RNA in most known examples.

\vskip .1in

\noindent Bondings are naturally partially ordered by inclusion, and the restriction to secondary structures on a fixed number of sites
thus has a geometric realization given by an arc complex for a polygon as was observed in [PW].
By the Classical Fact 11, this space of secondary structures is a sphere.  Statically, this gives combinatorial information
about the space of secondary structures, and dynamically, one imagines the evolution of a secondary structure in time as a dynamical
system on a suitable sphere.
Also proved in [PW] is the fact that an appropriate space of all RNA secondary structures is likewise homotopy equivalent to an arc
complex (for instance, the space of all RNA secondary structures with no non-trivial helices and a fixed number of unbonded
sites is again homeomorphic to the arc complex of a polygon), and a novel constructive proof (based on train tracks) of the Classical
Fact 11 is presented.
  
\vskip .1in

\noindent In this context, the planar sphericity conjecture has the immediate interpretation as capturing the
combinatorics of several interacting secondary structures on macromolecules in both the general and binary contexts.
Any model of full three-dimensional folding of systems of macromolecules  should be explicable in terms of arc families, and towards
this end, we next describe several recent models of ``pseudo-knotting'' of secondary structures. Among the many elaborations and
extensions of the basic thermodynamical model Mfold described before, [IS] and especially [XBTI] emphasize a related notion of
pseudo-knotting among other innovations.

\vskip .1in

\vskip .1in

\noindent The first model we propose is the three-dimensional
folding of a single binary macromolecule.

\vskip .09in

\leftskip .2in\rightskip .2in

\noindent {\bf Binary Chiral Model $F$}~Given a binary bonding ${\cal B}$ on $I_m$, construct a bordered surface
$F'\subseteq {\bf R}^3$ as follows. Attach to $I_m\times [0,1]\subseteq {\bf R}^2\subseteq {\bf R}^3$ a collection of bands in the
natural way, one band
$\beta _k$ for each bond $\{ i_i,j_k\} \in{\cal B}$, where the core of the band $\beta _k$ agrees with $C_{i_k,j_k}$; the
$z$-coordinates are adjusted in any convenient way to arrange that the bands are all disjointly embedded.  Consider a boundary
component $K$ of
$F'$; as illustrated in Figure 13, there are two cases for $K$: either the projection of $K$ bounds
a polygon, or it does not, and we 
adjoin to $F'$ an abstract polygon along each component of the former type to produce a surface $F$.
The boundary components of $F$ are
called ``pseudo-knots''. Notice that altering the lengths of helices or breaking helices by adding free sites leaves invariant the
topological type of
$F$, and $F$ is planar if and only if the folding is a secondary structure.

~~\vskip 1.3in

\hskip .3in{{{\epsffile{bnsffig1.epsf}}}}

\centerline{{\bf Figure 13}~~{Boundary components of $F'$.}}

\vskip .1in

\leftskip=0ex

\noindent The ``achiral model'' for binary folding allows each endpoint of a band to be attached either above the $x$-axis
or below the $x$-axis (so there are four ways to attach a band, and more data is required for each band).  In the same manner as
before, we may build a surface
$F'$ and cap off certain of its boundary components to produce another surface $F$, whose boundary components are the pseudo-knots. 
The chiral and achiral models for arbitrary macromolecules are defined similarly.   It is worth saying explicitly that there are
further elaborations of these four basic models (chiral/achiral for binary/arbitrary) tailored to other natural
biophysical circumstances.

\vskip .1in

\noindent We had already studied in [P2] a thermodynamic model with
Boltzmann weights given by simplicial coordinates; this physics project was undertaken as a pure mathematician
since the structure and attendant calculations were pleasing.  The achiral model in the binary case is closely
related to the thermodynamic model described in [P2].  (The two cyclic orderings on a set with three elements in [P2] corresponds
to the bond being attached above or below the $x$-axis.)

\vskip .1in

\vskip .1in
\noindent There is no doubt that the combinatorics of suitable models of macromolecules are described by arc complexes.  Whether
protein folding will be effectively predicted by {\sl ab initio} or by refined empirical motif recognition, a better understanding of
this combinatorics should be biologically useful. 

\vskip .1in

\noindent {\bf Speculative Remark }~As we have seen in $\S$2, the simplicial coordinates are effectively volumes in Minkowski space, so
in the spirit of [Kh] and [AS], the thermodynamic model in [P2] has ``the Boltzmann weights as volumes''.  Might the convex hull
construction itself give a method of predicting RNA or protein folding? 
The key question is how to
impose the primary structure of the macromolecule on the geometry.  Perhaps of utility in this regard, [KV] describes a hyperbolic
version of conformal field theory, which is especially interesting to us because the authors argue that it models suitable polymers;
furthermore
(without realizing it!), the authors employ
the special and unnatural case of decorated Teichm\"uller theory in which all of the horocycles have a common Euclidean radius
$\varepsilon$ in a model of upper half space, and then take the limit $\varepsilon\to 0$. 
Essentially the entire analysis of [KV] applies in our setting and suggests an explicit model of ``hyperbolic field theory'', where the
two-point function is the lambda length.

\vskip .09in

\vfill\eject

\noindent {\bf Bibliography}

\vskip .2in

\noindent [AS]~~M. Atiyah and P. Sutcliffe, ``The geometry of point particles''; hep-th/0105179.

\vskip .1in

\noindent [CS]   
M.\ Chas and D.\ Sullivan.
``String topology'', 
math.GT/9911159; to appear {\it Ann. Math}.

\vskip .1in

\noindent [CJ] R.\ L.\ Cohen and J.D.S.\ Jones,
{``A homotopy theoretic realization of string topology''},
Preprint math.GT/0107187.

\vskip .1in

\noindent [DM]~~P. Deligne and D. Mumford, ``Irreducibility of the space
of curves of a given genus'', {\it Publ. IHES} {\bf 36} (1979), 75-110.

\vskip .1in

\noindent [Ge]~ 
E.\ Getzler,
``Two-dimensional topological gravity and equivariant cohomology'', 
 {\it Comm.\ Math.\ Phys.\ } {\bf 163} (1994), 473--489.

\vskip .1in

\noindent [Ha] J. L. Harer, ``The virtual cohomological dimension of the mapping class group of an orientable surface, {\it Inv. Math.}
{\bf 84} (1986), 157-176.

\vskip .1in

\noindent [HM] J. H. Hubbard and H. Masur, ``Quadratic differentials and foliations'', {\it Acta Math.}{\bf 142} (1979), 221-274.

\vskip .1in

\noindent [K1] R. M. Kaufmann, ``On several varieties of cacti and their relations'', preprint MPI-2002-113, math.QA/0209131.

\vskip .1in

\noindent [K2]~---, this volume.

\vskip .1in

\noindent [KLP]~Ralph L. Kaufmann, Muriel Livernet, R. C. Penner, ``Arc operads and arc algebras'', 
preprint math.GT/0209132.

\vskip .1in

\noindent [Kh]~~A. Kholodenko, ``Kontsevich-Witten model from 2+1 gravity: new exact combinatorial solution'',
{\it Jour. Geom. and Phys.} {\bf 43} (2002), 45-91.

\vskip .1in

\noindent [Ko] M. Kontsevich, ``Intersection theory on the moduli space of curves and the matrix Airy function'', {\it Comm. Math.
Phys.} {\bf 147} (1992), 1-23.

\vskip .07in

\noindent [KV]~~~P. Kleban and I. Vassileva, ``Conformal field theory and hyperbolic geometry'', {\it Phys. Rev.
Lett.} {\bf 72} (1994) , 3929-3932.

\vskip .1in

\noindent [IS] H. Isambert and E. Siggia, {\it Proc. NAS USA} {\bf 97} (2002), 6515-6520.

\vskip .1in

\noindent [P1]~~ R. C. Penner, ``The decorated Teichm\"uller space of  punctured surfaces", 
{\it Comm. Math.  Phys.}  {\bf 113}   (1987),  299-339.

\vskip .1in

\noindent [P2]~~---, ``An arithmetic problem in surface geometry'', {\it
The Moduli Space of Curves}, Birk-h\"auser (1995), eds. R. Dijgraaf, 
C. Faber, G. van der Geer, 427-466.

\vskip .1in

\noindent  [P3]~---,``The simplicial compactification of Riemann's moduli space'', 
Proceedings of the 37th Taniguchi Symposium, World Scientific (1996), 237-252.

\vskip .1in

\noindent [P4]~---, ``Decorated Teichm\"uller theory of bordered surfaces'',
math.GT/0210326.

\vskip .1in

\noindent [P5]~---, ``Perturbative series and the moduli space of Riemann surfaces'' {\it Jour. Diff. Geom.}
{\bf 27} (1988), 35-53;``Weil-Petersson volumes'', {\it Jour. Diff. Geom.}
{\bf 35} (1992), 559-608; ``Universal constructions in Teichm\"uller theory'' {\it Adv. Math.} {\bf 98} (1993), 143-215; 
(with F. Malikov) ``The Lie algebra of homeomorphisms of the circle'', {\it Adv. Math.} {\bf 140} (1998), 282-322; 
``On Hilbert, Fourier, and wavelet transforms'', {\it Comm. Pure Appl. Math.} {\bf 55} (2002), 772-814. 

\vskip .1in

\noindent [PH]
R.\ C.\ Penner with J.\ L.\ Harer,
{\it Combinatorics of Train Tracks},
Annals of Mathematical Studies  125,  Princeton Univ. Press (1992, 2001).

\vskip .1in

\noindent [PW]~~R. C. Penner and M. S. Waterman  ``Spaces of RNA secondary structures",  
{\it Advances in 		Mathematics}  {\bf 101}  (1993), 31-49.

\vskip .1in

\noindent [St] K. Strebel, {Quadratic Differentials}, {\it Ergebnisse der Math.} {\bf 3:5}, Springer-Verlag, Heidelberg (1984).

\vskip .1in

\noindent [Su] D. Sullivan, ``Moduli of closed string operators and the homology action'' (2002); to appear in Abel Bicentennial
Proceedings.

\vskip .1in

\noindent [Vo] 
A.\ A.\  Voronov,
``Notes on universal algebra'',
math.QA/0111009.

\vskip .1in

\noindent [XBTI]  A. Xayaphoummine, T. Bucher, F. Thalmann, and H. Isambert, ``Prediction and statistics of pseudoknots in RNA
structures using exactly clustered stochastic simulations'', pre-print (2003).

\vskip .1in

\noindent [ZMT]~~M. Zuker, D. H. Mathews, and D. H. Turner, ``Algorithms and Thermodynamics for RNA secondary
structure prediction: a practical guide'', In {\it RNA Biochemistry and Biotechnology}, Kluwer (1999), 11-43.

\vskip .1in

\noindent [ZS]~~M. Zuker and D. Sankoff, ``RNA secondary structures and their prediction'', {\it Bull. Math.
Biol.} {\bf 46} (1984), 591-621.

\vskip .1in

\bye